\documentclass{amsart}

\usepackage{dsfont}
\usepackage[linesnumbered,boxed,vlined,ruled]{algorithm2e}
\usepackage{hyperref}
\usepackage{tikz}

\newcount\algocount
\def\looker{\ifx\next\end\else\ifx\next\lineending\vskip-4pt%
\indent\qquad\else\vskip0pt\global\advance\algocount by1\relax%
\indent\hskip12pt\hbox to0pt{\hss\the\algocount.}\quad\fi\fi}

\def\lineending{\parfillskip=0pt plus1fil\relax\futurelet\next\looker}

\def\thealgocount{\the\algocount}

{\obeylines\obeyspaces
\gdef\numberedalgorithm{%
\vskip\baselineskip\parindent=\thealgoindent\global\algocount=0\relax%
\def\note##1{\hfill##1\hskip-\parfillskip}\def\ {\qquad}%
\baselineskip=12pt\parskip=2pt\obeylines\obeyspaces%
\let^^M=\lineending\indent}}

\newif\ifcenteron
\newif\ifcodeon
\newif\ifneedspace
{\obeylines\obeyspaces

\gdef\codelooker{\ifx\next\codelineending\vskip-4pt%
\indent\else%
\ifx\next\end\vskip-4pt\let^^M\ \else%
\vskip0pt\indent\fi\fi}

\gdef\codelineending{\relax\futurelet\next\codelooker}

\gdef\algoindent#1{\gdef\thealgoindent{#1}}
\algoindent{0pt}

\gdef\algorithm{\medskip\global\codeontrue\parindent=\thealgoindent%
\let\ \qquad%
\def\note##1{\hfill##1\hskip-\parfillskip}%
\obeylines\obeyspaces\let^^M=\codelineending}

\gdef\spcodesamp{\baselineskip=12pt%
\parskip=2pt\obeylines\obeyspaces\let^^M=\codelineending}

\gdef\codebox#1{\let\go\relax%
\parindent=\thealgoindent%
\def\note##1{\hfill##1\hskip-\parfillskip}%
\ifcodeon\else\global\needspacetrue\let\go\spcodesamp\fi\go%
\global\dimen1=#1\relax\global\setbox0=\vbox\bgroup\indent}

\gdef\endcodebox{\ifneedspace\global\needspacefalse\vskip\baselineskip\fi%
\egroup%
\ifcenteron\vskip\baselineskip\hbox to\textwidth\bgroup\hss\vbox\bgroup\fi%
\vtop to \ht0{\vskip-3pt\hrule\hbox to\dimen1{\vrule height\ht0%
\hfill\vrule height\ht0 depth \dp0}\hrule\vss}%
\nobreak\vskip-\ht0\rlap{\hskip6pt\vbox{\unvbox0}}%
\ifcenteron\egroup\hss\egroup\global\centeronfalse\fi\algoindent{0pt}}

}

\def\change#1#2{\red{#1}\blue{#2}}
\def\change#1#2{#2}
\def\heron{(l_1+l_ 2+l_3)(l_1+l_2-l_3)(l_1-l_2+l_3)(-l_1+l_2+l_3)}

\title{Hyperbolic polyhedral surfaces with regular faces}
\author{Yohji Akama}
\address[Yohji Akama]{Mathematical Institute, Graduate School of Science, Tohoku University\\
6-3 Aoba, Sendai, 980-8578, Japan}
 \email{yoji.akama.e8@tohoku.ac.jp}

\author{Bobo Hua}
\address[Bobo Hua]{School of Mathematical Sciences, LMNS, Fudan University,
Shanghai, 200433, China}
\address{
Shanghai Center for Mathematical Sciences,
Fudan University, Shanghai,
200433,
China}
\email{bobohua@fudan.edu.cn}

\def\arccosh{\mathrm{arccosh}}
\def\implies{\;\Rightarrow\;}
\newcommand{\SP}{\mathds{S}}
\newcommand{\curvature}[1]{
 \ifnum #1=1
 5.74238\cdots\times 10^{-5}\else{%
 \ifnum #1=2
 1.83701\cdots\times 10^{-5}\else{%
 \ifnum #1=3
 7.85456\cdots\times 10^{-6}\else{%
 \ifnum #1=4
 3.49781\cdots \times10^{-6}\else{%
 \ifnum #1=5
 1.46541\cdots\times 10^{-6}\else{%
 \ifnum #1=6
 1.38186\cdots\times 10^{-8}\else{%
 \ifnum #1=7
 1.12992\cdots\times 10^{-8}\else{
 \ifnum #1=8
 4.70209\cdots\times 10^{-9} \else{
 \ifnum #1=9
 1.30776\cdots\times 10^{-9} \else{
 \ifnum #1=10
 4.96239\cdots\times 10^{-10} \else{
 \ifnum #1=11			% 7-7
 5.95826\cdots\times 10^{-2} \else{
 \ifnum #1=12                   % 8-9
 2.16500\cdots\times 10^{-3} \else{
 \ifnum #1=13			% 10-10
 1.04884\cdots\times 10^{-3}
\fi}
\fi}
\fi}
\fi}
\fi}
\fi}
\fi}
\fi}
\fi}
\fi}
\fi}
\fi}
\fi}
\def\aux{\alpha}

\def\lbFacialDeg{7}
\def\lbFacialDeg{16}
\def\ubFacialDeg{59}
\def\false{\textbf{false}}
\def\gonality{\B{2}{3}}
\def\true{\textbf{true}}

\newcommand\generalcurv[1]{
\ifnum #1>54
\curvature{10} \else{
\ifnum #1>44
\curvature{9} \else{
\ifnum #1>40
\curvature{8}  \else{
\ifnum #1>34
\curvature{7}\else{
\ifnum #1>20
\curvature{6} \else{
\ifnum #1>18
\curvature{5} \else{
\ifnum #1>16
\curvature{4} \else{
\ifnum #1>15
\curvature{3}
\fi}\fi}\fi}\fi}\fi}\fi}\fi}\fi}

\newcommand\generalp[1]{
\ifnum #1>54
(29,55,55) \else{
\ifnum #1>44
(35,38,43) \else{
\ifnum #1>40
(9,36,41)  \else{
\ifnum #1>34
(15,18,19) \else{
\ifnum #1>20
(16,16,21) \else{
\ifnum #1>18
(11,16,16) \else{
\ifnum #1>16
(11,13,15) \else{
\ifnum #1>15
(12,12,16) \else{
\ifnum #1>12
(8,10,13) \else{
\ifnum #1>10
(6,7,11) \else{
\ifnum #1>9
(5,6,10) \else{
\ifnum #1>7
(6,7,7) \else
(5,7,7)
\fi}\fi}\fi}\fi}\fi}\fi}\fi}\fi}\fi}\fi}\fi}\fi}

\newcommand\generalq[1]{
\ifnum #1>54
(32,43,55) \else{
\ifnum #1>44
(31,44,45) \else{
\ifnum #1>40
(10,23,26) \else{
\ifnum #1>34
(11,26,35) \else{
\ifnum #1>20
(14,20,20) \else{
\ifnum #1>18
(10,17,19) \else{
\ifnum #1>16
(11,12,17) \else{
\ifnum #1>15
(11,14,15) \else{
\ifnum #1>12
(7,13,13) \else{
\ifnum #1>10
(6,8,9) \else{
\ifnum #1>9
(5,7,8) \else{
\ifnum #1>7
(5,8,8)\else{
\ifnum #1>6
(6,6,7)\fi}\fi}
\fi}\fi}\fi}\fi}\fi}\fi}\fi}\fi}\fi}\fi}\fi}

\def\B#1#2{B_{#1}(#2)}

\def\gap{\kappa}
\def\gp#1{\epsilon_{\max}^{#1}}
\def\neararea#1{\wt{\area^{#1}_{\max}}}
\def\EqRef#1{(\ref{#1})}

\newcommand{\Cat}{{\mathrm{CAT_{{loc}}}(-1)}}
\newcommand{\SH}{S^{\mathds{H}^2}}

\newcommand{\pt}{\mathrm{Pttn}}

\newcommand{\R}{{\mathds R}}

\newcommand{\area}{\mathrm{Area}}

\newcommand{\cri}{\mathrm{cri}}
\newcommand{\HP}{\mathds{H}^2}
\newcommand{\wt}{\widetilde}

\newcommand{\PC}{\mathcal{PC}_{>0}}

\newcommand{\ST}{\mathcal{T}_{S_g}}
\def\adp#1{\mathcal{A}{\mathrm{d}}\mathcal{P}_{#1}}
\newcommand{\ad}{\mathcal{A}{\mathrm{d}}\mathcal{P}}
\def\adp#1{\mathcal{A}{\mathrm{d}}\mathcal{P}}

\newcommand{\Ng}{\mathcal{NC}_{<0}}
\newcommand{\NC}{\mathcal{NC}_{\leq -1}}
\newcommand{\NCg}{\mathcal{NC}_{\leq -1}^{g}}
\newtheorem{theorem}{Theorem}
\newtheorem{corollary}[theorem]{Corollary}
\newtheorem{df}[theorem]{Definition}
\newtheorem{eg}[theorem]{Example}
\newtheorem{lemma}[theorem]{Lemma}
\newtheorem{problem}[theorem]{Problem}
\newtheorem{proposition}[theorem]{Proposition}
\begin{document}

\begin{abstract}
 We study hyperbolic polyhedral surfaces with faces
 isometric to regular hyperbolic polygons satisfying that the total angles at vertices are at least $2\pi.$
The combinatorial information of these surfaces is shown to be
 identified with that of Euclidean polyhedral surfaces with negative
 combinatorial curvature everywhere.  We prove that there is a gap
 between areas of non-smooth hyperbolic polyhedral surfaces and the area
 of smooth hyperbolic surfaces. The numerical result for the gap is
 obtained for hyperbolic polyhedral surfaces, homeomorphic to the double
 torus, whose 1-skeletons are cubic graphs.
\end{abstract}
\thanks{The first author is partially supported by JSPS KAKENHI [grant Number JP16K05247].}
\maketitle
 \section{Introduction}\label{sec:intro}

The combinatorial curvature for planar graphs, as the generalization of the Gaussian curvature for surfaces was introduced by \cite{Nevanlinna70,Stone76,Gromov87,Ishida90}. Many interesting geometric and combinatorial results have been obtained since then, see e.g.
\cite{Zuk97,Woess98,Higuchi01,BP01,HJL02,LPZ02,
HigS03,SunYu04,RBK05,BP06,DeVosMohar07,
ChenChen08,Zhang08,Chenbf09,Keller10,KP11,Keller11,Oh17,Gh17}.

Let $(V,E)$ be an undirected, locally finite, simple graph with the set of vertices $V$ and the set of edges $E.$ The graph $(V,E)$ is called \emph{semiplanar} if it is topologically embedded into the surface $S,$ see \cite{HJL15}. We write $G=(V,E,F)$ for the combinatorial structure, or the cell complex, induced by the embedding where $F$ is the set of faces, i.e. connected components of the complement of the embedding image of the graph $(V,E)$ in the target. Two elements in $V,E,F$ are called \emph{incident} if the closures of their images have non-empty intersection. We say that a graph $G$ is a \emph{tessellation} of $S$ if the following hold, see e.g. \cite{Keller11}:
\begin{enumerate}
\item \label{i}Every face is homeomorphic to a disk whose boundary consists of finitely many edges of the graph.
\item \label{ii} Every edge is contained in exactly two different faces.
\item \label{iii} For all two faces whose closures have non-empty intersection, the intersection is either a vertex or an edge.
\end{enumerate} In this paper, we only consider tessellations and call them semiplanar graphs for the sake of simplicity. For each semiplanar graph, we always assume that for each vertex $x$ and face $\sigma,$ $$\deg(x)\geq 3,\ \mathrm{deg}(\sigma)\geq 3$$ where $\deg(\cdot)$ denotes the degree of a vertex or a face.
For each semiplanar graph $G$, the \emph{combinatorial curvature} at the
vertex is defined as
\begin{equation}\label{def:comb}\Phi(x)=1-\frac{\deg(x)}{2}+\sum_{\sigma\in
F \colon x\in \overline{\sigma}}\frac{1}{\deg(\sigma)},\quad (x\in V),\end{equation} where the summation is taken over all faces $\sigma$ incident to $x.$  To digest the definition, we endow the ambient space of $G$ with a canonical piecewise flat metric and call it the (regular) \emph{Euclidean polyhedral surface}, denoted by $S(G)$: Replace each face by a regular Euclidean polygon of same facial degree and of side length one, glue them together along their common edges, and define the metric on the ambient space via gluing metrics, see \cite[Chapter~3]{BuragoBuragoIvanov01}. It is well-known that the generalized Gaussian curvature on a Euclidean polyhedral surface, as a measure, concentrates on the vertices. And one is ready to see that the combinatorial curvature at a vertex is in fact the mass of the generalized Gaussian curvature at that vertex up to the normalization $2\pi,$ see e.g. \cite{MR2127379,HJL15}. For each finite semiplanar graph embedded into a surface $S,$ the Gauss-Bonnet theorem, see e.g. \cite[Theorem~1.2]{DeVosMohar07}, reads as
\begin{equation}\label{eq:gbf}
\sum_{x\in V}\Phi(x)=\chi(S),
\end{equation} where $\chi(\cdot)$ is the Euler characteristic of the surface.

We denote by $$\Ng:=\{G=(V,E,F)\colon \Phi(x)<0,\forall x\in
V\}$$  the class of semiplanar graphs with negative combinatorial curvature everywhere. There are many examples in the class $\Ng.$

We review some known results on the class $\Ng.$ For each semiplanar graph $G$ in $S=\R^2,$ there hold various isoperimetric inequalities, see e.g. \cite{Zuk97,Woess98,Higuchi01,HJL02,LPZ02,HigS03,KP11,Keller11}.
The following proposition is proved by Higuchi.
\begin{proposition}[\cite{Higuchi01}, Proposition~2.1]\label{prop:Higuchi} \begin{equation}\label{eq:1806}x\in V, \Phi(x)<0\implies \Phi(x)\leq-\frac{1}{1806}.\end{equation}
Equality holds if and only if $\deg x=3$ and $x$ is incident to  3-, 7-, and 43-gons.\end{proposition}

\section{Hyperbolic polyhedral surfaces and main results}\label{s:hyperbolic setting}

In this paper, we study hyperbolic polyhedral surfaces with faces isometric to regular hyperbolic polygons. Let $\HP$ be the simply connected hyperbolic surface of constant curvature $-1.$ We only consider regular hyperbolic polygons in the hyperbolic space $\HP$ with at least three sides. For each $n\geq 3$ and $a>0,$ there is a  regular hyperbolic $n$-gon in $\HP$ of side length $a,$ denoted by $\Delta_n(a),$ which is unique up to the hyperbolic isometry. Analogous to Euclidean polyhedral surfaces, we define hyperbolic polyhedral surfaces associated to a semiplanar graph. For each semiplanar graph $G=(V,E,F)$ and $a>0,$ we replace each face by a regular hyperbolic polygon in $\HP$ of side length $a,$ and glue them together along their common edges. This induces a metric structure on the ambient space of $G,$ called \emph{hyperbolic polyhedral surface} of $G$ with side length $a$ and denoted by $\SH_a(G).$ In the literature, hyperbolic polyhedral surfaces are well studied for circle packings on triangulations of surfaces whose side lengths are induced by the radii of circles, which indicates that the triangles are not necessarily regular {in that setting}, see e.g. \cite{Thurston76,ChowLuo03,LuoGuDai08,GeXu16DGA,GeXu17IMRN}. In our setting, we consider general polyhedral surfaces with arbitrary faces, but restrict to those with constant side length everywhere, focusing on the combinatorics of such surfaces.

We denote by $\area_a(G)$ the area of $\SH_a(G),$ which is possibly infinite.  For each $x\in V,$ the total angle at $x$ measured in $\SH_a(G)$ is denoted by $\theta_a(x),$ and the angle defect at $x$ is defined as
$$K_a(x)=2\pi-\theta_a(x),\qquad a>0.$$ By the hyperbolic geometry, we yield the Gauss-Bonnet theorem on hyperbolic polyhedral surfaces. {Since the argument is standard, we omit the proof, see e.g. \cite{AkamaHuaSu18}.}
\begin{theorem}\label{thm:GaussBonnet} For each finite semiplanar graph $G=(V,E,F)$ embedded into a surface $S$ and $a>0,$
\begin{equation}\label{GaussBonnet}-\area_a(G)+\sum_{x\in V}K_a(x)=2\pi \chi(S).\end{equation}
\end{theorem}

We say that a geodesic metric space $(X,d)$ locally has the (sectional) curvature bounded above by $-1$ in the sense of Alexandrov if it satisfies the Toponogov triangle comparison property with respect to the hyperbolic space {$\HP$}, and denote by $\Cat$ the set of such spaces, see e.g. \cite{BuragoBuragoIvanov01}.  It is well-known that $\SH_a(G)\in \Cat$ if and only if
$$K_a(x)\leq 0,\quad \forall x\in V.$$ We denote by $$\NC:=\{G=(V,E,F)\colon  \SH_a(G)\in \Cat\ \mathrm{\ for\ some\ } a>0 \}$$ the class of  semiplanar graphs whose hyperbolic polyhedral surface of certain side length has the curvature locally bounded above by $-1$ in the Alexandrov sense. We will prove that for each finite semiplanar graph $G$,
$$G\in\NC\iff G\in \Ng,$$ see Proposition~\ref{prop:pcspc}. This suggests a possible way to study the class $\NC$ by known results on the class $\Ng,$ where the latter refers to the Euclidean setting.

From now on, we only consider finite semiplanar graphs.

\begin{df}\label{def:critical side length} Let $G=(V,E,F)\in\NC$.
\begin{enumerate}
\item We define the pattern of $x\in V$ by  $$\pt(x):=(\mathrm{deg}(\sigma_1),\mathrm{deg}(\sigma_2),\cdots,\mathrm{deg}(\sigma_{N})),$$ where $\{\sigma_i\}_{i=1}^N$ are the faces incident to $x$ ordered by $\mathrm{deg}(\sigma_1)\leq\mathrm{deg}(\sigma_2)\leq\cdots\leq\mathrm{deg}(\sigma_{N}),$ and $N=\deg(x).$

\item We define \emph{the critical side length of  a vertex $x\in V$},  \emph{the critical  side length of  $G$},  and \emph{the critical  area of  $G$}, respectively by
\[
 \begin{array}{rl}
 a_c(x)&:=\max\left\{a>0\colon K_a(x)\leq 0\right\},\\
 a_c(G)&:=\min_{x\in V}a_c(x),\ \mbox{and}\\
 \area^{\cri}(G)&:=\area_{a_c(G)}(G).
 \end{array}
\]
 \end{enumerate}
\end{df}

One is ready to see that $\area_a(G)$ is monotonely increasing in $a,$ which implies that $$\area^{\cri}(G)=\max\{\area_a(G)\colon\SH_a(G)\in \Cat\}.$$

Let $G$ be a semiplanar graph embedded into $S_g,$ the closed orientable surface of genus $g\geq 2.$ {In particular, for $g=2,$ $S_2$ is called the double torus.}
We say that a semiplanar graph $G$ \emph{admits a hyperbolic tiling with
regular hyperbolic polygons} if there is a hyperbolic tiling with
regular hyperbolic polygons of $S_g$ equipped with a smooth hyperbolic metric, whose semiplanar graph structure is isomorphic to $G,$ and denote by $\ST$ the set of such semiplanar graphs. One is ready to prove the following proposition.

\begin{proposition}\label{prop:tiling} For each semiplanar graph $G$ embedded in $S_g$, the following are equivalent:
\begin{enumerate}
\item $G\in \ST.$
\item There is some $a>0$ such that $\SH_a(G)$ is isometric to $S_g$ equipped with a smooth hyperbolic metric.
\item $a_c(x)=a_c(y)$ for all $x,y\in V.$
\item $\area^{\cri}(G)=2\pi(2g-2).$
\end{enumerate}
\end{proposition}

For $g\geq 2,$ we denote by $$\NCg:=\{G=(V,E,F)\in \NC\colon S(G)=S_g,\,
\# V<\infty\}$$ {the set of} finite semiplanar graphs embedded into
$S_g$ whose hyperbolic polyhedral surface $\SH_a(G)\in \Cat,$ for some
$a>0.$

\begin{eg}\label{eg:Bolza}\rm The \emph{Bolza
	   surface}~\cite{Bolza1887} is a compact Riemann surface of
	   genus~2. 

 \begin{enumerate}
  \item The Bolza surface is a hyperbolic surface that can be defined by
	a subgroup of the $(2,3,8)$ triangle
	group~\cite[Section~3]{KKSV16}~(See also \cite[Lemma~2.2]{KS06}).
	That is, the hyperbolic triangle of inner angles $\pi/2,\pi/3$, and $\pi/8$
	tiles the Bolza surface by reflection of the triangle on the
	edges. The set of the centers of the incircles of the
	triangles induces a tessellation
	$G_B$ of $S_2$. Actually, $G_B\in\mathcal{T}_{S_2}$, because the side length of the tessellation
	is the diameter $a$ of the incircles and $K_v(a)=0$
	for each vertex $v$ of $G_B$.
	
  \item  We consider the
 fundamental domain of the Bolza surface in the Poincar\'e disk where the
 opposite sides of the octagon are
 identified in the octagon in Figure~\ref{fig:Bolza}.
       Figure~\ref{fig:Bolza} is a Delaunay triangulation of the Bolza
	 surface where the points of the same name are identified.
	 Lemma~\ref{lem:A23} proves that the Voronoi
	 tessellation corresponding to the Delaunay triangulation is in
	$\NC^2\setminus\mathcal{T}_{S_2}$.        Figure~\ref{fig:Bolza} designates
the degrees of the 15 vertices of the triangulation, and enumerates the
	34 faces of the triangulation.

Figure~\ref{fig:Bolza} is drawn by a program the first author wrote with ``2D Periodic Hyperbolic Triangulations'' package~\cite{iordanov17:_implem_delaun_trian_bolza_surfac} of
 \texttt{CGAL~5.4}. 
 \end{enumerate}
\def\ps{{\large a}}
\def\pr{{\large b}}
\def\pu{{\large c}}
\def\pv{{\large i}}
\def\pw{{\large e}}
\def\px{{\large f}}
\def\py{{\large g}}
\def\pz{{\large h}}
\def\pp{{\large d}}
 \begin{figure}[ht]\centering
       \begin{tikzpicture}
	\begin{scope}[scale=.6]
	 \node (russel) at (-6,0)  {
	 %.55*.6
	 \includegraphics[scale=.33]{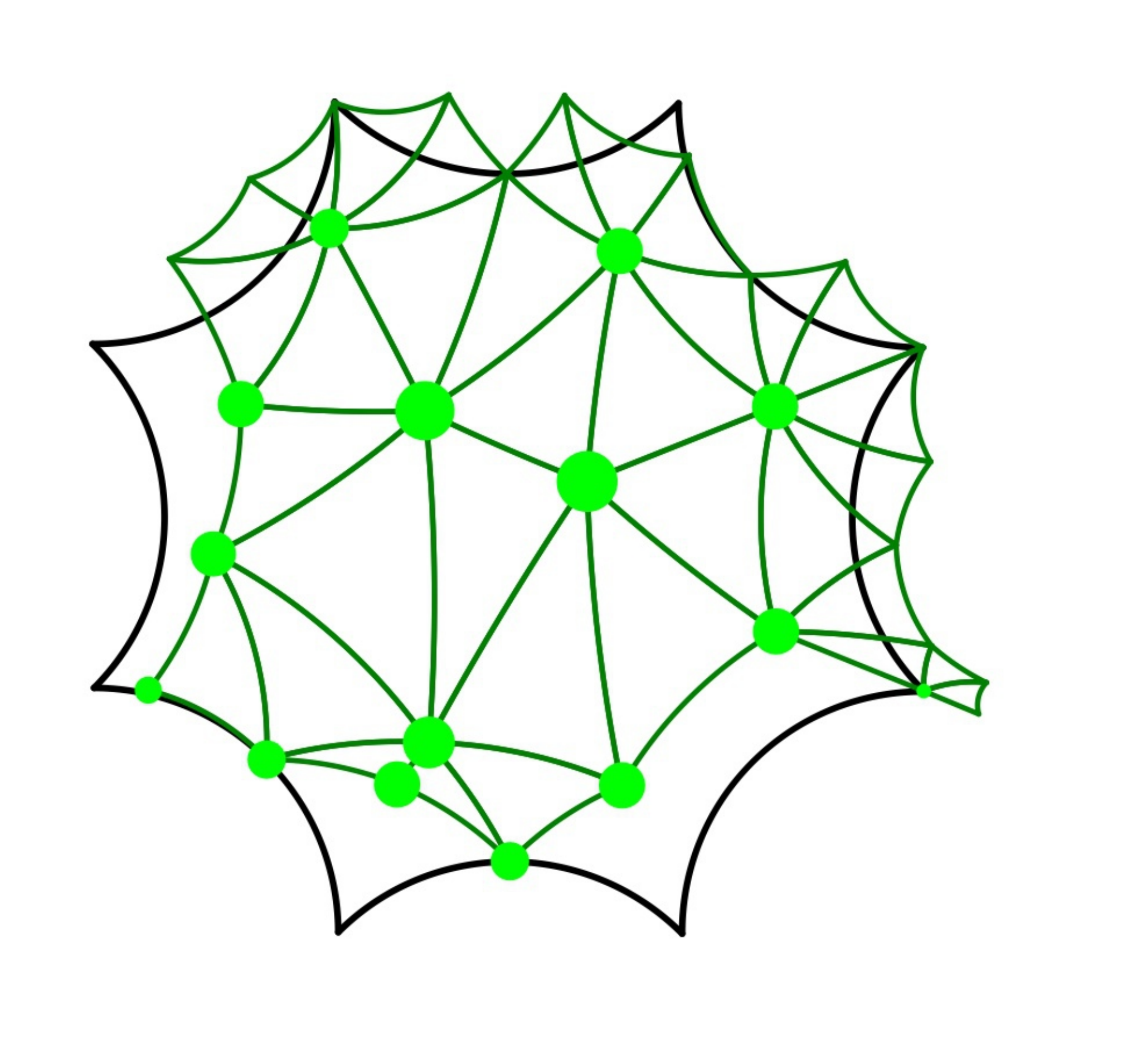}}; %.55*.6

	 %
	 % The north arc
       \node at (-9.9,7.5) {\pr}; %vertex b outside the fundamental domain
       \node at (-8,7.6) {\pu};   %vertex c outside the fundamental domain
       \node at (-7,6.5) {\pp};   %vertex d outside the fundamental domain
       \node at (-5.9,7.6) {\pw}; %vertex e outside the fundamental domain

	 % The north eastern arch
       \node at (-3.4,6.1) {\px}; %vertex f outside the fundamental domain
       \node at (-2.5,4.8) {\py}; %vertex g outside the fundamental domain
       \node at (-1,4.8) {\pu};   %vertex c outside the fundamental domain
       \node at (.3,3.4) {\pr};   %vertex b outside the fundamental domain

	 % The eastern arch
       \node at (.5,1.2) {\pv}; %vertex i outside the fundamental domain
       \node at (0,-.1) {\pz}; %vertex h outside the fundamental domain

       \node at (.5,-1.6) {\px}; %vertex f outside the fundamental domain
       \node at (1.5,-2.4) {\pw}; %vertex e outside the fundamental domain
       \node at (1.3,-3.2) {\pv};  %vertex i outside the fundamental domain
       \node at (0,-4) {\pr};   %vertex b inside the fundamental domain

       \node at (-2.5,-2.3) {\ps}; %vertex a inside the fundamental domain
       \node at (-4.5,-4.9) {\pw}; %vertex e inside the fundamental domain

       \node at (-6.9,-6.3) {\pp};   %vertex d inside the fundamental domain
       \node at (-8.9,-4.8) {\pu};   %vertex c intside the fundamental domain

       \node at (-11.2,-4.4) {\py}; %vertex g inside the fundamental domain
       \node at (-12.9,-3.3) {\px}; %vertex f inside the fundamental domain

       \node at (-12.3,.3) {\pz}; %vertex h inside the fundamental domain
       \node at (-12,1.9) {\pv}; %vertex i inside the fundamental domain

	 % North western arc
       \node at (-12.7,4.9) {\pw}; %vertex e outside the fundamental domain
       \node at (-11.3,6.4) {\ps}; %vertex a outside the fundamental domain
       
       \node at (-9.9,5) {\Huge 7};
       \node at (-5,4.7) {\Huge 7};
       
       \node at (-11.3,2.2) {\Huge 6};
       \node at (-8.3,2) {\Huge 7};
       \node at (-2.4,2.1) {\Huge 8};

       \node at (-5.5,.9) {\Huge 6};

       \node at (-11.8,-.3) {\Huge 7};

       \node at (-2.4,-1.7) {\Huge 7};

       \node at (-12.8,-2.6) {\Huge 6};
       \node at (0,-3.3) {\Huge 7};

       \node at (-10.9,-3.7) {\Huge 6};
       \node at (-8.1,-3.5) {\Huge 7};
       \node at (-5,-4.2) {\Huge 9};

       \node at (-8.7,-4.2) {\Huge 6};
       
       \node at (-6.9,-5.5) {\Huge 6};

       \node at (-11.3, 5) {5};
       \node at (-10.4, 6) {8};
       \node at (-9.4, 6.2) {10};
       \node at (-8.5, 5.5) {34};
       \node at (-6, 5.5) {3};
       \node at (-4.9, 5.6) {4};
       \node at (-4,5) {33};

       \node at (-11, 3.3) {16};
       \node at (-9.5, 3.3) {11};
       \node at (-8.5, 4.5) {12};
       \node at (-6.5, 4.5) {15};
       \node at (-3.3,3.5) {7};
       \node at (-2.3,3.4) {13};

       \node at (-1.3,3) {18};
       \node at (-1.3,2) {17};

       \node at (-10,1.6) {19};
       
       \node at (-6.5, 2.5) {29};
       \node at (-4.5, 2.5) {6};
       \node at (-1.4,1.1) {32};

       \node at (-9.3,-.3) {26};
       \node at (-7.3,-.3) {30};
       \node at (-3.8, .5) {1};
       \node at (-1.8, 0) {25};
       
       \node at (-12,-2.3) {21};
       \node at (-10,-2.3) {20};
       \node at (-6.5,-2.3) {28};
       \node at (-4.5,-.9) {23};
       \node at (-1.2,-1.1) {2};
       \node at (-.8,-2) {27};
       \node at (-.8,-2) {27};
       
%       \node at (.47,-2.3) {9}; %face 9 outside the fundamental domain
%       \node at (.49,-2.7) {14}; %face 14 outside the fundamental domain

       \node at (-9.3,-3.8) {24};
       \node at (-8,-4.3) {31};
       \node at (-6.5,-4.3) {22};
	\end{scope}
       \end{tikzpicture}
	    \caption{\label{fig:Bolza} The faces 9 and 14 are the
  triangles bef and bei, respectively.}
   \end{figure}
\end{eg}

We are interested in critical areas of semiplanar graphs in $\NCg$. In particular, we propose the following problem.
\begin{problem}\label{prob:minarea} What is the constant
$$\wt{\area_{\max}}:=\sup_{G\in \NCg\setminus \ST} (\area^{\cri}(G))\ ?$$
\end{problem}

By Proposition~\ref{prop:tiling}, $$\area^{\cri}(G)=2\pi(2g-2),\quad \mathrm{for\ each\ } G\in \ST.$$ One is ready to see that for each $G\in\change{\NC}{\NCg}$ by Gauss-Bonnet formula \EqRef{GaussBonnet},
$$\area^{\cri}(G)\leq 2\pi(2g-2).$$ Hence graphs in the class $\ST$ attain the maximal critical area in the class $\change{\NC}{\NCg}.$ The problem is to determine the maximum of critical areas in the class $\change{\NC}{\NCg}$ except hyperbolic tilings $\ST.$
We prove the following gap for the areas.
\begin{theorem}\label{thm:eps} For each $g\geq 2,$ there is a constant $\epsilon>0,$ depending on $g,$ such that
$$\wt{\area_{\max}}\leq 2\pi(2g-2)-\epsilon.$$
\end{theorem}
\begin{proof}
We will show that there are only finitely many graphs in $\NCg,$ see Corollary~\ref{coro:finite} in Section~\ref{s:generalsurface}. Hence the gap for the areas exists by the assertion $(4)$ in Proposition~\ref{prop:tiling}.
\end{proof}

Concerning with the result in Theorem~\ref{thm:eps}, we aim to obtain a quantitative estimate of $\wt{\area_{\max}}$ in the following.

 \begin{df}\label{def:tentative gap}Let
  \begin{eqnarray*}
   A_g(N):=\{G=(V,E,F)\in  \NCg\setminus \ST\colon \deg(x)\le N, \forall x\in V\ \},\\
  \neararea{N}:=\sup_{G\in A_g(N)} \area^{\cri}(G)\qquad(3\le
   N\le\infty).
  \end{eqnarray*}
 \end{df}  
\change{}{
Lemma~\ref{lem:A23} proves that Figure~\ref{fig:Bolza} is indeed a member of $A_2(3)$.
 In Corollary~\ref{coro:finite}, we establish
$\#\left(\NCg\right)<\infty$ for each $g\geq2.$
} Note that
\begin{enumerate}
\item  $\neararea{N}$ is nondecreasing in $N$, and
\item  $\neararea{\infty}=\wt{\area_{\max{}}}$.
\end{enumerate}

\begin{theorem}\label{thm:gap}{The following is the
  computation result when we represent numbers with 40
  decimal digit in the computation:} For $g=2$,
 \[  \neararea{3}\le 4\pi-  {\generalcurv{\ubFacialDeg}}.
  \]
\end{theorem}

Let $g=2.$ Our strategy is to reformulate the problem to a new one
which can be estimated by local arguments. We rewrite
\begin{equation}\label{eq:a11}\wt{\area_{\max}}=4\pi-\gp{N},\end{equation} where
{$$\gp{N}:=\inf_{G\in A_2(N)} \left(4\pi - \area^{\cri}(G)\right).$$}

We  call this \emph{the gap} between the
maximal critical area $4\pi$ and other critical areas. By the
Gauss-Bonnet theorem, \EqRef{GaussBonnet}, we have
\begin{equation}
  \gp{N}=\inf_{G\in A_2(N)} \sum_{\mbox{$x\in V$}}(-K_{a_c(G)}(x)).\label{total angle defect}
\end{equation}
Hence for the upper bound
estimate in Theorem~\ref{thm:gap}, it suffices to obtain the lower bound estimate of the absolute value of total angle
defect for these semiplanar graphs. This new problem fits to local
arguments, and we prove the results by enumerating all cases, see Section~\ref{s:double torus} and Appendix.

{The paper is organized as follows: In next section, we introduce basic
properties of hyperbolic polyhedral surfaces and prove
Proposition~\ref{prop:tiling}. We bound combinatorial quantities of
semiplanar graphs in the class $\NCg,$ which are embedded in general
surfaces, in Section~\ref{s:generalsurface}. In Section~\ref{s:double
torus}, we refine the above estimates for cubic graphs embedded in the
double torus and prove Theorem~\ref{thm:gap}. In Section~\ref{sec:future
work}, we propose some further works. We present the algorithm for Theorem~\ref{thm:gap} in Appendix.}

\section{Preliminaries}
Let $G=(V,E,F)$ be a  semiplanar graph. Two vertices are called \emph{neighbors} if there is an edge connecting them. We denote by $\deg(x)$ the degree of a vertex $x,$ i.e. the number of neighbors of a vertex $x,$ and by $\deg(\sigma)$ the degree of a face $\sigma,$ i.e. the number of edges incident to a face $\sigma$ (equivalently, the number of vertices incident to $\sigma$).

For each regular hyperbolic $n$-gon of side length $a>0,$ $\Delta_n(a)$, in $\HP,$ we denote by $\beta=\beta_{n,a}$ the inner angle at corners of the $n$-gon. By the hyperbolic geometry,
\[\cosh\frac{a}{2}\sin\frac{\beta}{2}=\cos\frac{\pi}{n}.\] One is ready to see that $\beta_{n,a}$ is monotonely decreasing in $a$ and
\begin{equation}\label{eq:theta}\lim_{a\to 0}\beta_{n,a}=\frac{n-2}{n}\pi,\quad \lim_{a\to \infty}\beta_{n,a}=0\end{equation} where the right side of the first equation is the inner angle of a regular $n$-gon in the plane.
The area of $\Delta_n(a)$ is given by
$$\area(\Delta_n(a))=(n-2)\pi-n\beta.$$

\begin{proposition}\label{prop:pcspc} For each finite semiplanar graph $G$,
$$G\in\NC\iff G\in \Ng,$$
\end{proposition}
\begin{proof} For each $G\in \NC,$ there is $a>0$ such that $\SH_a(G)\in \Cat.$ For each face $\sigma\in F$ with $\deg(\sigma)=n,$ we know that the inner angle of $\Delta_n(a)$ is less than that of a Euclidean $n$-gon. Consider the Euclidean polyhedral surface $S(G).$ Note that the total angle at each vertex in $S(G)$ is greater than that in $\SH_a(G).$ This yields that $G\in \Ng.$ Thus $\NC\subset \Ng.$

For the other direction, let $G\in \Ng.$ Note that by \EqRef{eq:theta}, for each vertex $x\in V,$ there is a small constant $a(x)$ such that the total angle \change{$\theta_{a(x)}(x)<2\pi.$}{$\theta_{a(x)}(x)>2\pi.$} Since the graph is finite, we can choose a small constant $a$ such that
$$\SH_a(G)\in \Cat,$$ which proves that $G\in \NC.$ This proves the proposition.
\end{proof}

For each $a>0,$ the total angle at the vertex $x$ in $\SH_a(G)$ for some
graph $G$ is given by $$\theta_a(x)=\sum_{i=1}^N2\arcsin
\frac{\cos{\frac{\pi}{f_i}}}{\cosh{\frac{a}{2}}},\qquad\left((f_1,\ldots,f_N)=\pt(x)\right).$$
As $a>0$, $\cosh (a/2)>1$. Thus $\arcsin(\cos(\pi/f_i)/\cosh(a/2))$ is
defined.  To determine the critical side length of the vertex, $a_c(x)$
is the unique solution to the following equation $$\theta_a(x)=2\pi.$$

Now we prove Proposition~\ref{prop:tiling}.
\begin{proof}[Proof of Proposition~\ref{prop:tiling}] $(1)\Longleftrightarrow (2):$ This is trivial.

$(2)\Longrightarrow (3):$ This follows from the monotonicity of $\theta_a(x)$ in $a$ for each $x\in V.$

$(3)\Longrightarrow (4):$ For $S_{a_c(G)}(G),$ it is smooth at each vertex, hence it is locally isometric to a domain in $\HP.$ This implies that $S_{a_c(G)}(G)$ is isometric to a hyperbolic surface. By the Gauss-Bonnet formula \EqRef{GaussBonnet}, $\area^{\cri}(G)=2\pi(2g-2).$

$(4)\Longrightarrow (2):$ We know that the area of $S_{a_c(G)}(G)$ is $2\pi(2g-2).$ By the Gauss-Bonnet formula \EqRef{GaussBonnet}, $K_{a_c(G)}(x)=0$ for all $x\in V.$ Hence $S_{a_c(G)}(G)$ is a smooth hyperbolic surface.
\end{proof}

\begin{df}\label{def:monotonic} Let
  $p=(f_1,\ldots,f_N)$ and $q=(g_1,\ldots,g_M)$ be  two nondecreasing integer sequences with $f_i,g_j\ge3$.
  	 \begin{eqnarray*}
	 \Phi(p)&:=&1 -  \sum_{i=1}^N\left(\frac{1}{2} - \frac{1}{f_i}\right).\\
	   K_a(p)&:=&2\pi -
	 \sum_{i=1}^N	 2\arcsin\frac{\cos\frac{\pi}{f_i}}{\cosh\frac{a}{2}}\qquad(a>0).\\
	 a_c(p)&:=&\max\{a>0 \colon K_a(p)\le 0\}.\\
	 p\preceq q &:\iff& N\leq M,\
 \mathrm{and}\ f_{N-i}\leq g_{M-i}\ ( 0\leq\forall i\leq N-1).
\end{eqnarray*}
 Let $\prec$ be the strict part $\preceq\setminus =$.
\end{df}

\begin{lemma}\label{lem:monotonic}Let
  $p,q$ be as in Definition~\ref{def:monotonic}.
\begin{enumerate}
    \item \label{assert:curvature defined} $K_a(p)$ is
	 a strictly increasing, continuous function of $a$.

	 \item \label{assert:critical side length}Suppose $\Phi(p)<0$. Then
 $a_c(p)$ is well-defined, and is the unique solution
	       $a$ such that $K_a(p)=0$ and $a>0$.

 \item \label{assert:ac} Suppose that $3\le f_1\le f_2\le
       f_3$ are  nondecreasing integers,  $l_i=\cos(\pi/f_i)$
       $(i=1,2,3)$, and $H=\heron$. Then
   \begin{eqnarray}
a_c(f_1,f_2,f_3)
  &=
\arccosh \left(\frac{8\prod_{i=1}^3
l_i^2 }{H
	}
  -1 \right).
       \label{def:alpha}
\end{eqnarray}

$\cosh \left(a_c(f_1,f_2,f_3)/2\right)$ is the half of the volume of Euclidean cuboid
 of sides $l_i$ $(i=1,2,3)$ divided by the  area of Euclidean
 triangle  of sides $l_i$ $(i=1,2,3)$.

 \item \label{assert:ac monotone}$p \change{\preceq}{\prec} q\implies  a_c(p)\change{\leq}{<}
	      a_c(q).$
 \item \label{assert:Phi monotone} $f_i\le f_i'$ $(1\le
		    i\le N)\implies \Phi(f_1,\ldots,f_N)\ge \Phi(f'_1,\ldots,f'_N)$.
  \end{enumerate}
 \end{lemma}
  \begin{proof}
   The assertion~\EqRef{assert:curvature defined} is clear.

   The assertion~\EqRef{assert:critical side length} follows
from the assertion~\EqRef{assert:curvature defined} and  $\lim_{a\to\infty} K_a(p)=2\pi$ and $\lim_{a\to0}K_a(p)=2\pi
   \Phi(p)<0$.

The assertion~\EqRef{assert:ac}. The denominator in the argument of $\arccosh$ is
   well-defined, because
   \begin{equation}
\frac{1}{2}=
  \cos\frac{\pi}{3}\le \cos\frac{\pi}{f_i}<1.    \label{withintwice}
\end{equation}
   The equation~\EqRef{def:alpha} is proved as follows: By
   Lemma~\ref{lem:monotonic}~\EqRef{assert:critical side length},
   $a=a_c(p)$ satisfies $2\pi = \sum_{i=1}^3 2\arcsin\left(
   {\cos(\pi/f_i)}\,/\,{\cosh(a/2)} \right)$. Hence,
 \[
\cos\left(\pi - \arcsin \frac{l_1}{\cosh\frac{a}{2}}\right) =
\cos\left( \arcsin \frac{l_2}{\cosh\frac{a}{2}} +\arcsin \frac{l_3}{\cosh\frac{a}{2}}\right) .
 \]
By the addition formula for cosines,
 \[
  -\sqrt{1-\frac{l_1^2}{\cosh^2\frac{a}{2}}}
 =
\sqrt{\left(1-\frac{l_2^2}{\cosh^2\frac{a}{2}}\right)
 \left(1-\frac{l_3^2}{\cosh^2\frac{a}{2}}\right)}
 -\frac{l_2 l_3}{\cosh^2\frac{a}{2}}.
 \]
 By multiplying both sides by $\cosh^2(a/2)$,
 \[
  -\cosh \frac{a}{2} \sqrt {\cosh^2\frac{a}{2} - l_1^2}
 =
 -l_2  l_3
 +\sqrt{\left(\cosh^2\frac{a}{2}  - l_2^2\right)
        \left(\cosh^2 \frac{a}{2} - l_3^2\right) }.
 \]
 By squaring both sides and then moving terms,
\[      
  -2 l_3^2 l_2^2+(l_3^2+l_2^2-l_1^2)\cosh^2\frac{a}{2}
 =
-2 l_2 l_3 \sqrt{\left(l_2^2-\cosh^2\frac{a}{2}\right)\left( l_3^2-\cosh^2\frac{a}{2}\right)}.     \]
By subtracting the square of the right side from the square of the
 left side,
 \[
- H
  \cosh^4\frac{a}{2} +4   l_1^2 l_2^2 l_3^2
\cosh^2\frac{a}{2} = 0.
\]
 Because $a$ is a real number, $\cosh\frac{a}{2}>0$. Thus,
$\cosh^2\frac{a}{2} = {4 \prod_{i=1}^3 l_i^2}/{H}$.
 From $\cosh^2(a/2)=(1+\cosh a)/2$, we obtain the desired equation.  By Heron's area formula for Euclidean triangles, we
   have the last sentence.

   The assertion~\EqRef{assert:ac monotone}.  For $p \change{\prec}{\preceq} q$,
 $K_a(p)\change{\geq}{>} K_a(q)$ by Definition~\ref{def:monotonic}.
 In particular, $K_{a_c(p)}(p)>K_{a_c(p)}(q)$.
 By the second assertion of this Lemma, $0=K_{a_c(q)}(q)=K_{a_c(p)}(p)$.
 Thus, $K_{a_c(q)}(q)>K_{a_c(p)}(q)$.
  By the first assertion of this Lemma, $a_c(q)>a_c(p)$.

The assertion~\EqRef{assert:Phi monotone} is obvious.
  \end{proof}

\section{Number of vertices, vertex degrees and face degrees}\label{s:generalsurface}
Let $S=S_{g}$ be a closed orientable surface of genus $g\geq 0.$
Let $G$ be a tessellation of $S.$  For each $G\in \NC,$ by the
Gauss-Bonnet formula~(Theorem~\ref{thm:GaussBonnet}), $g\geq 2.$ Then Gauss-Bonnet formula
\EqRef{eq:gbf}, the result of Higuchi, and Proposition~\ref{prop:Higuchi} yield the following result.
\begin{proposition}\label{prop:gest} For $G\in \NCg,$
$$\# V\leq 3612(g-1).$$
\end{proposition}
\begin{proof} By \EqRef{eq:gbf} and \EqRef{eq:1806}, $$\# V \frac{1}{1806}\leq \sum_{x\in V}|\Phi(x)|\leq 2g-2.$$ This yields the result.
\end{proof}

By this proposition, we have the following corollary.
\begin{corollary}\label{coro:finite} For each $g\geq 2,$  $\#\left(\NCg\right)<\infty.$
\end{corollary}

Moreover, we can estimate the vertex degree.
\begin{proposition}\label{prop:vertexdeg} Let $G\in \NCg$\change{}{ with $g\ge2$}.  For each $x\in V,$
$$\deg(x)\leq 12 g-7.$$
\end{proposition}
\begin{proof} Suppose it is not true. Let $x_0$ be a vertex such that $\deg(x_0)\geq12 g-6.$ Then
 $$\Phi(x_0)\leq 1 - \frac{\deg{(x_0)}}{2} + \sum_{x_0\in \overline{\sigma},
 \sigma\in F} \frac{1}{\deg(\sigma)} \leq 1 - \frac{\deg{(x_0)}}{2} +
\frac{ \deg{(x_0)}}{3},$$ {since $\deg(\sigma) \ge 3$ $(\forall \sigma\in F)$}. Thus, $\Phi(x_0)\le
 1-\frac{\deg(x_0)}{6}\le  2-2g.$ Since there are at least two vertices in the graph,  by Proposition~\ref{prop:pcspc} the Gauss-Bonnet formula \EqRef{eq:gbf} yields a contradiction to
 $G\in \NCg.$
\end{proof}

\medskip
We recall the result of \cite[Table~1]{Higuchi01}.
 \begin{lemma}[\protect{\cite[Table~1]{Higuchi01}}] \label{Hig} For
  $f_1,\ldots,f_N$~$(3\le f_1\le f_2\le\cdots\le f_N)$,
  $\Phi(f_1,\ldots,f_N)<0$  if and only if
  \begin{itemize}
   \item $N\ge7$;
   \item $N=6$ and $(f_1,\ldots,f_N)\succeq (3,3,3,3,3, 4)$;
   \item $N=5$ and $(f_1,\ldots,f_N)\succeq p$ for some $p=(3,3,3,3,7)$,
	 $(3, 3, 3, 4, 5)$,  $(3, 3, 4, 4, 4)$;
   \item $N=4$ and $(f_1,\ldots,f_N)\succeq p$ for some $p=(3,3,4,13)$,
	 $(3,3,5,8)$, $(3,3,6,7)$, $(3,4,4,7)$, $(3,4,5,5)$,
	 $(4,4,4,5)$; or
   \item $N=3$ and
 $(f_1,f_2,f_3) \succeq  p$ for some $p=
(3, 7, 43)$,  $(3, 8, 25)$, $(3, 9, 19)$, $ (3, 10, 16)$, $(3, 11, 14)$, $ (3, 12, 13)$,
	 $(4, 5, 21)$, $(4, 6, 13)$, $(4, 7, 10)$,
	 $(4, 8, 9)$, $(5, 5, 11)$, $ (5, 6, 8)$,
	 $(5, 7, 7),$ {or} $(6, 6, 7)$.
	 \end{itemize}
\end{lemma}

\begin{proposition}\label{prop:facedegree} For $G\in \NCg$\change{}{ with $g\ge2$},
$$\deg(\sigma)\leq 84g-43,\quad \forall \sigma\in F.$$
\end{proposition}
\begin{proof} Let $\sigma$ be the face with maximal facial
 degree {$f=\deg(\sigma).$}

On the one hand, by the tessellation properties \EqRef{i} and \EqRef{ii} in the
 introduction, $\# F\geq 2$ and there is at least one vertex $y$ which
 is not on the boundary of $\sigma.$ By Proposition~\ref{prop:pcspc},
 $\Phi(y)<0$.
Thus, by the Gauss-Bonnet formula \EqRef{eq:gbf},
\begin{equation}\label{eq:125eq1:21-3-19}\sum_{x\in
\overline{\sigma}}\Phi(x)>2-2g.\end{equation}
On the other hand,  by Table~1 in \cite{Higuchi01}, i.e.,
 Lemma~\ref{Hig}, we have:
  \begin{equation}  %
\pt(x)\succeq   (\overbrace{3,\ldots,3}^{\deg(x)-1},\ f),\ 
   (3,3,4,f),\
  (3,7,f), \ \mathrm{or}\ (4,5,f),\quad (x\in\overline{\sigma}). \label{ptx:21-3-19}
  \end{equation}
 Indeed, for $N:=\deg(x)=3$, $\pt(x)\succeq(3,7,f)$ for
 $\pt(x)=(3,\ldots)$, and $\pt(x)\succeq(4,5,f)$ otherwise. For $N=4$,
 $\pt(x)\succeq(3,3,4,f)$, and for $N\ge5$,
 $\pt(x)\succeq(3,\ldots,3,f)$ as in \EqRef{ptx:21-3-19}.
 
 The maximum
 $\Phi(p)$ where $p$
 ranges over the four tuples in the right side of \EqRef{ptx:21-3-19}, is $\Phi(3,7,f)$.  By \EqRef{ptx:21-3-19} and Lemma~\ref{lem:monotonic}~\EqRef{assert:Phi
 monotone}, 
 \[\Phi(x)\le \Phi(3,7,f)=-\frac{1}{42}+\frac{1}{f},\quad
 (x\in\overline{\sigma}).\] By summing both  sides  over the $f$
 vertices $x\in\overline{\sigma}$, 
\[\sum_{x\in \overline{\sigma}}\Phi(x)\leq
-\frac{f}{42}+1.\]  Thus, the conclusion follows from \EqRef{eq:125eq1:21-3-19}.
\end{proof}

\section{Cubic graphs embedded in the double torus}\label{s:double torus}
In this section, we study {cubic graphs, i.e. regular graphs of vertex
degree $3$, embedded into the double torus $S_2.$}

\begin{lemma}\label{lem:A23}
 The Delaunay triangulation  of Bolza surface (Figure~\ref{fig:Bolza})
 induces a Voronoi tessellation which is a member of $A_2(3)$.
\end{lemma}
 \begin{proof}We can check that
  the vertex patterns of the Voronoi tessellation are as in
  Table~\ref{tbl:A23}. Here, the numbering of the vertices in Table~\ref{tbl:A23} corresponds
  to the numbering of the faces of the Delaunay triangulation in Figure~\ref{fig:Bolza}.
\begin{table}[ht]\centering
 \begin{tabular}{c|ccc||c|ccc||c|ccc||c|ccc}
 no.& & & & no.& & & & no.& & & &  no.& & &  \\
  \hline
  1 & 6 & 7 & 8	& 2 & 6 & 7 & 7	& 3 & 6 & 7 & 9	& 4 & 6 & 7 & 9\\
5 & 7 & 7 & 9	& 6 & 6 & 7 & 8	& 7 & 6 & 7 & 8	& 8 & 7 & 7 & 7\\
9 & 6 & 7 & 9	& 10 & 6 & 7 & 7	& 11 & 6 & 7 & 7	& 12 & 6 & 7 & 7\\
13 & 6 & 6 & 8	& 14 & 6 & 7 & 9	& 15 & 6 & 7 & 7	& 16 & 6 & 7 & 9\\
17 & 6 & 7 & 8	& 18 & 6 & 7 & 8	& 19 & 6 & 7 & 7	& 20 & 6 & 7 & 7\\
21 & 6 & 6 & 7	& 22 & 6 & 7 & 9	& 23 & 6 & 7 & 9	& 24 & 6 & 6 & 7\\
25 & 7 & 7 & 8	& 26 & 7 & 7 & 7	& 27 & 6 & 7 & 7	& 28 & 6 & 7 & 9\\
29 & 6 & 7 & 7	& 30 & 6 & 7 & 7	& 31 & 6 & 6 & 7	& 32 & 6 & 7 & 8\\
33 & 6 & 6 & 7	& 34 & 6 & 6 & 7	
\end{tabular}\caption{\label{tbl:A23}}\end{table}
All the vertices have negative combinatorial curvature, so the graph is in $\Ng^2$ by Lemma~\ref{Hig}, and thus in $\NC^2$.
The vertex pattern $(6,7,8)$ of vertex no.~1 is strictly greater than
 the vertex pattern $(6,7,7)$ of the vertex no.~2, so the two vertices have different critical side lengths by
 Lemma~\ref{lem:monotonic}~\EqRef{assert:ac monotone}.
 Therefore, the graph is not in $\mathcal{T}_{S_2}$ by Proposition~\ref{prop:tiling}.
 \end{proof}

       \begin{df}\label{def:B}For $g\ge2$,
	     \[
	      \B{g}{3}:=\max\left\{ \deg\sigma\colon (V,E,F) \in \NCg,\ \deg x=3\;(\forall x \in V),\ \sigma\in F\right\}.
	     \]
      \end{df}

\begin{proposition}\label{prop:Bg3 at least 7}
 $\B{g}{3}\ge7$.
\end{proposition}
\begin{proof}
 Owing to Proposition~\ref{prop:pcspc}, $\NC^2\subseteq\Ng$. By Lemma~\ref{Hig},
 every vertex pattern $(f_1, f_2, f_3)$ of $G\in \NC^2\subseteq\Ng$ is greater than or equal to one of the 14 triples $(3,7,43),\ldots,(6,6,7)$.  Hence, $f_3\ge43,\ldots,$ or $7.$ Therefore, the maximum facial degree $\B{g}{3}$ of $G$ is the maximum of such $f_3$, and is greater than or equal to $7$.
\end{proof}
By Example~\ref{eg:Bolza}~(1), the tessellation
$G_{B}\in\mathcal{T}_{S_2}$ of the Bolza surface provides a lower bound 16 on $\B{2}{3}$. Moreover, for cubic graphs,
we can improve the upper bound $84g - 43$ on the
facial degree that Proposition~\ref{prop:facedegree} provides, to $40g - 21$.
\begin{proposition}\label{prop:59}For $g\ge2$,
 $\B{g}{3}\le 40g - 21$.
 \end{proposition}
 \begin{proof} Suppose that $G=(V,E,F)\in \NCg$ and
  $\deg x=3$ for all $x\in V$.
Suppose that $\sigma$ is a face of maximal facial degree $f$.
By Lemma~\ref{Hig}, we know that for each $x\in \overline{\sigma},$ which
					    is not of pattern $(3,7,f)$
		or $(3,8,f),$ \begin{equation}\label{eq:p13}\Phi(x)\leq
			      -\frac{1}{20}+\frac{1}{f}.\end{equation}
  Set
  \[
  W:=\{y\in \overline{\sigma}\colon \pt(y)=(3,7,f)\ \mathrm{or}\
  (3,8,f)\}.\]
  Let $y\in W.$ We denote by $w$ the vertex satisfying $w\in \overline{\sigma},$ $w\sim y$ and $w,y$ are incident to a common triangle, see Figure~\ref{fig378}.
 \begin{figure}[htbp]
   \begin{center}
     \begin{tikzpicture}
      \draw[thick]
      (5.4,5)--(6.8,5)--(7.6,6.1)--(7.2,7.4)--(6.2,8.1)--(4.9,7.4)--(4.8,6)--cycle;

      \draw[thick] (4.8,6)--(4.2,5);

      \node at (4.8,5.4) {$3$};
      \node at (4.2,4.7) {$w$};
      \node at (5.4,4.7) {$y$};
      \node at (4.6,6)   {$z$};
      \node at (4.8,3.4) {$f$};

      \node at (6,6.5) {$7,8$};
      
      \draw[thick]
      (5.4,5)--(6.8,5)--(7.6,4.6)--(8.2,3.8)--(8.4,2.8)--(7.4,1.8)--(5.5,1.8);

      \draw[thick,dashed] (5.5,1.8)--(4.5,1.8);

      \draw[thick]
      (4.5,1.8)--(2.8,1.8)--(1.9,2.5)--(2,3.8)--(2.8,4.8)--(4.2,5)--(5.4,5);

      \draw[thick]
      (4.9,7.4)--(4.3,8.1)--(3.3,8.4)--(2.2,8);

      \draw[dashed,thick] (2.2,8)--(1.6,6.9);

      \draw[thick]
      (1.6,6.9)--(1.7,5.6)--(2.8,4.8);
     \end{tikzpicture}
  \caption{
  $\pt(y)=(3,7,f)$ or $(3,8,f).$}
\label{fig378}
   \end{center}
\end{figure} We define a map, \begin{eqnarray*}T:&&W\to V\cap \overline{\sigma},\\
&&y\mapsto w.\end{eqnarray*} We denote by $z$ the vertex adjacent to $w$
  and $y,$ by $\tau$ the face incident to $w,$ which is \change{not
  $\sigma$ and}{neither $\sigma$ nor}
  the triangle $\Delta_{y,z,w}.$ Since $\Phi(z)<0,$ by Lemma~\ref{Hig}
\begin{equation}
  \deg(\tau)\geq 25.\label{ineq:25}\end{equation} Therefore, by $\pt(w)=(\change{\ldots}{3},\deg\tau,f)$,
	     $w\notin W$. Hence \begin{equation}TW\cap W=\emptyset.\label{disjoint}\end{equation}
  Let $u$ be the vertex on $\overline{\sigma}$
  adjacent to $w,$ which is not $y.$  By \EqRef{ineq:25}, $\pt(u)=(\ldots,\deg \tau, f)$ is {neither}
  $(3,7,f)$ {nor}
  $(3,8,f).$  Thus $u\not\in W.$ Hence $T^{-1}(w)=y$ and the map $T$ is injective. By \EqRef{ineq:25},
\begin{equation}\label{eq:p12}\Phi(Ty)=\Phi(w)=-\frac12+\frac13+\frac{1}{\deg(\tau)}+\frac1f\le
 -\frac{19}{150}+\frac1f\quad(y\in
W).\end{equation}
  Obviously
  \begin{equation}
   \Phi(y)\le
   -\frac{1}{2}+\frac{1}{3}+\frac{1}{7}+\frac{1}{f}\quad(y\in W). \label{ineq:y}
  \end{equation}
By \EqRef{disjoint},
\[\sum_{x\in \overline{\sigma}}\Phi(x)=\sum_{y\in
 W}(\Phi(y)+\Phi(Ty))+\sum_{
x\in \overline{\sigma}\setminus(W\cup TW)}\Phi(x).\]
  Hence by {\EqRef{ineq:y}, \EqRef{eq:p12}, and \EqRef{eq:p13}},
  $\sum_{x\in \overline{\sigma}}\Phi(x)$ is at most
  \begin{eqnarray*}
&&(\# W)\left(-\frac{1}{42}+\frac{1}{f}\right) +
 (\# (T W))\left(-\frac{19}{150}+\frac{1}{f}\right) \\
   &&+(f - 2(\# W)) \left(-\frac{1}{20}+ \frac{1}{f}\right)
   \leq-{\frac {53\,(\# W)}{1050}}-\frac{f}{20}+1,
   \end{eqnarray*}
by \EqRef{disjoint} and by
$
  \# W=\# (TW)$ that follows from the injectivity of $T.$
In addition, by the tessellation properties \EqRef{i} and \EqRef{ii} in the
 introduction, there is at least one vertex which
 is not on the boundary of $\sigma,$ as in the proof of Proposition~\ref{prop:facedegree}. This yields that
  $$\sum_{x\in V}\Phi(x)<\sum_{x\in \overline{\sigma}}\Phi(x).$$

  Therefore
  by the Gauss-Bonnet formula \EqRef{eq:gbf},
  \begin{equation}
   2-2g=\sum_{x\in V}\Phi(x)<\sum_{x\in \overline{\sigma}}\Phi(x)\leq
   -{\frac {53\,(\# W)}{1050}}-\frac{f}{20}+1. \label{pr}
  \end{equation} By $\# W\ge0$, \EqRef{pr}, and $f\in\mathds{Z}$, we
  have  $f\le 40g - 21.$
 \end{proof}

\change{There is no classification of $\NC^2$,
 and $\B{2}{3}$ has not been determined yet.}{}
\begin{df} \label{df:ad}A nondecreasing integer sequence $p=(f_1,f_2,f_3)$ with $f_i\ge3$ is called \emph{admissible}, if
$\Phi(p)<0$ and $f_3\le \B{2}{3}$.  The set of admissible $p$
 is  denoted by $\ad$.
\end{df}
\change{}{By the finiteness of the set $\ad$, we can define the following:}
\begin{df}\label{def:k}
\[
\gap:=  \max\left\{K_{a_c(p)}(q)\ \colon
  K_{a_c(p)}(q)<0,\  p\in\ad,\ q\in \ad,\ a_c(p)<a_c(q) \right\}.
 \]
\end{df}

Recall from Section~\ref{s:hyperbolic setting}
\[
 \gp{N}=\inf_{G\in A_2(N)} \left(4\pi - \area^{\cri}(G)\right)=\inf_{G\in A_2(N)} \sum_{x\in V}(-K_{a_c(G)}(x)).
 \]

\begin{lemma}\label{lem:main}For $g\ge2$, $\gp{3}\ge-\gap$.\end{lemma}

 \begin{proof} By Corollary~\ref{coro:finite}, there are only finitely many $G$ in {$A_g(3)$}.
  Suppose $G=(V,E,F)\in A_g(3)$ and $a_c(G)=a_c(x)$ with $x\in V$. Then there is $y\in V$ such that $a_c(y)>a_c(x)$, by Proposition~\ref{prop:tiling} and Definition~\ref{def:critical side length}.
 By Lemma~\ref{lem:monotonic}, $K_{a_c(G)}(y)<K_{a_c(G)}(x)=0$. Thus,
\begin{eqnarray}
 \gp{3}&\ge & \inf\left\{  - K_{a_c(G)}(y)\colon {G\in A_g(3)},\ y\in
		   V,\ a_c(y)>a_c(x) \right\}\nonumber\\
& =&-\max\left\{   K_{a_c(G)}(y)\colon {G\in A_g(3)},\ y\in V,\ a_c(y)>a_c(x) \right\}.\nonumber\end{eqnarray}
 Suppose that $y\in V$ attains the maximum.  By ${G\in A_g(3)}$ and $a_c(x)<a_c(y)$,
 $G\in\NCg$.
 By Proposition~\ref{prop:pcspc}, $G\in\Ng$ and thus $\Phi(x), \Phi(y)<0$.  Therefore both of $\pt(x)$ and $\pt(y)$ are in $\ad$.
  Hence, the conclusion follows from Definition~\ref{def:k}.
 \end{proof}

\begin{proof}[Proof of Theorem~\ref{thm:gap}] By \EqRef{eq:a11}, \EqRef{total
  angle defect} and Lemma~\ref{lem:main}, it suffices to calculate
 $\gap.$
By Example~\ref{eg:Bolza}~(1) and Proposition~\ref{prop:59},
 \change{we have $\B{2}{3}\leq \ubFacialDeg$}{the constant $\B{2}{3}$ is
 in an interval $[\lbFacialDeg,\,\ubFacialDeg]$.  If $\B{g}{3}$ is
 estimated large, so are $\ad$ and $\kappa$ by Definition~\ref{df:ad} and Definition~\ref{def:k}}. For each
 $\B{2}{3}\in[\lbFacialDeg,\, \ubFacialDeg]$, we computed $\gap$ by
 using Maple with the algorithm given in Appendix. The result is
 Table~\ref{tbl:gap}. The largest $\gap$ is
 $$\gap=K_{a_c(\generalp{\ubFacialDeg})}(\generalq{\ubFacialDeg})=-\generalcurv{\ubFacialDeg}.$$
 This proves the theorem.
\end{proof}

 \begin{table}[htbp]\begin{center}
 \begin{tabular*}{\textwidth}{@{\extracolsep{\fill}}cccc}
  Range of $\B{2}{3}$&$\gap=K(a_c(p),\,q)$&$p$&$q$ \cr
    \hline
  16--16&$-\curvature{3}$&$\generalp{16}$&$\generalq{16}$\cr
  17--18&$-\curvature{4}$&$\generalp{17}$&$\generalq{17}$\cr
  19--20&$-\curvature{5}$&$\generalp{19}$&$\generalq{19}$\cr
  21--34&$-\curvature{6}$&$\generalp{21}$&$\generalq{21}$\cr
  35--40&$-\curvature{7}$&$\generalp{35}$&$\generalq{35}$\cr
  41--44&$-\curvature{8}$&$\generalp{41}$&$\generalq{41}$\cr
  45--54&$-\curvature{9}$&$\generalp{45}$&$\generalq{45}$\cr
  55--59&$-\curvature{10}$&$\generalp{55}$&$\generalq{55}$\cr
 \end{tabular*}

 \caption{ \label{tbl:gap}}		     \end{center}
\end{table}

\section{Future work}\label{sec:future work}

We recall the results in \cite{AkamaHuaSu18} for regular spherical polyhedral surfaces.
For each finite planar graph, it associates with some metric spaces, called
\emph{regular spherical polyhedral surfaces}, by replacing faces with
regular spherical polygons in the unit sphere and gluing them
edge-to-edge. We consider the class of planar graphs which admit
spherical polyhedral surfaces with the curvature bounded below by 1 in
the sense of Alexandrov, i.e.\@ the total angle at each vertex is at most
$2\pi$. We classify all spherical tilings with regular spherical
polygons, i.e.\@ total angles at vertices are exactly $2\pi$. We prove
that for each graph in this class which does not admit a spherical
tiling, the area of the associated spherical polyhedral surface with the
curvature bounded below by 1 is at most $4\pi-1.6471\cdots\times
10^{-5}$.  In other words,  $\wt{\area_{\max}}\le
4\pi-1.6471\cdots\times 10^{-5}$ for regular spherical polyhedral
surfaces.

In the spherical case, there are vertex types $p$ such that the
combinatorial curvature $\Phi(p)$ of $p$ is positive and
$K(a_c(p),p)>0$.
In the hyperbolic case,  $K(a_c(p),p)=0$ for every
vertex type $p$ with $\Phi(p)<0$.
We conjecture that $\wt{\area_{\max}}$ for regular
hyperbolic polyhedral surfaces with $g=2$ is greater than
$\wt{\area_{\max}}$ for regular spherical polyhedral surfaces, and is
close to $4\pi$.

For $$\PC:=\{G=(V,E,F)\colon \Phi(x)>0,\forall x\in V\},$$
DeVos and Mohar~\cite{DeVosMohar07} proved that any graph $G\in \PC$ is
finite, which solves a conjecture of Higuchi~\cite{Higuchi01}, see
\cite{Stone76,SunYu04} for early results. It would be interesting to classify $\PC$.  For the set
\[
 P:=\{\; G\; \colon\;  G\in \PC\mbox{ is neither a prism nor an antiprism}\;\},
\]
DeVos and Mohar proved that $\# P<\infty$ and asked the number
\[
C_{\SP^2}:=\max_{ (V,E,F) \in P}\# V.
\]
 For the lower bound estimate of $C_{\SP^2}$, large
examples in this class are constructed \cite{RBK05,MR2836763,Gh17,Ol17}, and finally some examples possessing $208$
vertices were found.  DeVos and Mohar~\cite{DeVosMohar07} showed
that $C_{\SP^2}\leq 3444,$ which was improved to $C_{\SP^2}\leq 380$ by
Oh~\cite{Oh17}. By a refined argument, in~\cite{Gh17}, Ghidelli
completely solved the problem.
  \begin{theorem}\label{thm:208}
   \begin{enumerate}
\item \label{assert:208} $C_{\SP^2}=208.$
      \cite{Gh17,MR2836763}

  \item
\label{thm:41upper}
$\max\{ \deg(\sigma)\colon \sigma\in F, (V,E,F)\in P\}\leq 41.$  \cite{Gh17}
\end{enumerate}\end{theorem}
Comparing with these results, we propose to determine
$$\max_{(V,E,F)\in \NCg} \# V \quad \mathrm{and}\quad \B{g}{3}\ \mathrm{for}\ g\ge2.$$

\change{a lower bound of $B_2(3)$, we want to find a tessellation
  is provided by a polyhedral embedding of the
snark $sn7$ in $S_2$,   in the proof of
Proposition~\ref{prop:11}.}{}

\change{We hope these results are useful to give a lower bound of $B_g(3)$.}{}

 Hyperbolic polyhedral surfaces enjoy more combinatorial
structures than the spherical ones by Proposition~\ref{prop:pcspc}. It
is a challenge to compute effective numerical results for ``almost''
hyperbolic tilings, i.e., non-smooth hyperbolic polyhedral surfaces.

\section*{Declaration of competing interest}
The authors declare that they have no known competing financial
interests or personal relationships that could have appeared to influence
the work reported in this paper.

\section*{Acknowledgements}
 The second author thanks Feng Luo for many helpful suggestions on
 hyperbolic polyhedral surfaces.
 We thank the anonymous referee for valuable comments and suggestions to
 improve the writing of the paper.

\appendix

\section{The algorithm}

We will try to numerically compute \change{$\kappa$.}{the following
   constant of Definition~\ref{def:k}:
   \[
\gap=  \max\left\{K_{a_c(p)}(q)\ \colon\ (p,q)\in
					   \ad^2,\
  K_{a_c(p)}(q)<0,\  \ a_c(p)<a_c(q) \right\}.
   \]}

   We note two difficulties:
    \begin{enumerate}
         \item It is expensive to check the conditions
	       for all the $(p,q)\in \ad^2$.

    \item
Since $|\gap|$ would be minute,
    the numerical computation of $\gap$ would suffer from
    \emph{loss of significance}~(so-called \emph{catastrophic
    cancellation})~\cite{Knuth:1997:ACP:270146}.

 \emph{We do not know whether the
	algorithm with 40 decimal digit
	representation of numbers is enough to confirm the value of
	$\gap$.} If the true value of $\gap$ has
	modulus less than $10^{-40}$, the
	algorithm misses the true value of
	$\gap$, because it computes $\gap$ as $2\pi$
	minus a total angle $\vartheta$ where the decimal
	representations of $2\pi$ and $\vartheta$ coincide up to the
	first 40 digits. By the same reasoning, if the true value of
	$\gap$ is indeed {$K_{a_c\generalp{\ubFacialDeg}}\generalq{\ubFacialDeg}$, then the number of
	significant digits of $\gap=-\generalcurv{\ubFacialDeg}$ is at most $40-\change{11}{10}$.}
    \end{enumerate}

To go around the first difficulty, i.e., to save the computation time, we first note that $\gap$ is
       the maximum of the following two:
	       \begin{eqnarray}
%		&&
		 \max\{
					   K_{a_c(p)}(q)\ \colon\ (p,q)\in
					   \ad^2,
					   K_{a_c(p)}(q)<0,\ a_c(p)<a_c(q),\ p\prec		q\}.\quad\label{sss}		\\
%		&&
		 \max\{
					   K_{a_c(p)}(q)\ \colon\ (p,q)\in
					   \ad^2,
					   K_{a_c(p)}(q)<0,\ a_c(p)<a_c(q),\ p\not\prec q\}.\quad\label{A.1.5}
\end{eqnarray}
	       Then,
	       \begin{itemize}
		\item
	       to totally dispense with checking many $(p,q)$ for
			       \EqRef{sss}, we will prove that the
			       following attains \EqRef{sss}:
	       \begin{equation}
\label{attain}					       (\tilde p,\ \tilde
				       q):=\left((B_2(3)-1, B_2(3),
			       B_2(3)),\ (B_2(3), B_2(3),
	       B_2(3))\right)\in\ad^2.\end{equation}

		\item
	       To ease the computation for \EqRef{A.1.5}, the following inexpensive
		      constraint is a necessary condition of
	       $K_{a_c(p)}(q)<0$ where $p=(f_1,f_2,f_3)$ and $q=(g_1,g_2,g_3)$:

		      \begin{df}\label{def:cncv}Let ${cncv(f_1,f_2,f_3,g_1,g_2,g_3)}$ be
\begin{eqnarray*}
      &\neg(f_1\ge g_1 \ \mbox{and}\  f_2-g_2{\ge}g_3-f_3  >0 ) \\
 \mbox{and}\  &\neg(f_2\ge g_2 \ \mbox{and}\  f_1-g_1{\ge}g_3-f_3  >0 ) \\
 \mbox{and}\  &\neg(f_3\ge g_3 \ \mbox{and}\  f_1-g_1{\ge}g_2-f_2  >0 ).
\end{eqnarray*}
		      \end{df}
	       \end{itemize}

The proofs of Lemma~\ref{lem:beta prop} and Lemma~\ref{aux} both depend on mean value theorem for
the following function:
	\begin{df}\label{df:beta}
	 \[
	     \beta(x,\,a):=2\arcsin
  \frac{\cos\frac{\pi}{x}}{\cosh\frac{a}{2}} \quad(3\le x\le\gonality, a>0).
	 \]
	\end{df}
	  For each integer $x\ge3$,  $\beta(x,\,a)$ is the inner angle $\beta_{x,a}$ of
  a regular hyperbolic $x$-gon of side length $a$.

\begin{lemma}\label{lem:beta prop}\label{lem:concave}\change{}{Let $3\le x\le \B{2}{3}$ and
 $a>0$. Then}
  \begin{enumerate}\item \label{assert:increasing and concave in x}$\beta(x,\,a)$ is increasing and concave in $x$.
   \item \label{assert:monotone in a} $\frac{\partial\beta}{\partial x}(x,a)$ is positive and is decreasing in $a$.
	 \end{enumerate}
\end{lemma}
 \begin{proof}
$\frac{\partial\beta}{\partial x}$ is $2\,{\pi \sin \left( {\frac {\pi }{x}} \right) {{x^{-2}
  u(x,a)^{-1/2}}}}>0$ where $$u(x,a):=- \left( \cos  \frac {\pi }{x}  \right) ^{2
}+ \left( \cosh \frac{a}{2}  \right) ^2>0.$$ As $u(x,a)$ is increasing
  in $a$,
$\frac{\partial\beta}{\partial x}$  is  decreasing in $a$. This
  establishes the second assertion.
 $\frac{\partial^2\beta}{\partial x^2}(x,a)$ is the product of a
  positive number $ u(x,a) ^{-3/2}{x}^{-4}>0$ and
\begin{eqnarray*}
-4\,\sin \left( {\frac {\pi }{x}} \right)  u(x,a)
 \pi \,x
 -2\,\cos \left( {\frac {\pi }{x}} \right) {\pi }^2 \left(
(\cosh \frac{a}{2})^2 -1 \right)<0 .\end{eqnarray*} Thus, $\frac{\partial^2\beta}{\partial x^2}(x,a)<0$\change{ for
  $x\ge3$}{}. \change{}{This establishes the first assertion and completes the proof of Lemma~\ref{lem:concave}.}
\end{proof}

\medskip
\change{}{By Definition~\ref{def:monotonic},  $K_{a_c(p)}(q)$ is
represented in the following form, because
Lemma~\ref{lem:monotonic}~\EqRef{assert:critical side length} implies
 $K_{a_c(p)}(p) = 2\pi - \sum_{i=1}^N \beta(f_i,a_c(p)) = 0$.}
\begin{lemma}\label{aux}
For $p=(f_1,\ldots,f_N)$, $q=(g_1,\ldots,g_M)\in\ad$,
   \[
    K_{a_c(p)}(q) =  \sum_{i=1}^N \beta(f_i,a_c(p)) -
    \sum_{i=1}^M \beta(g_i,a_c(p)).
   \]
\end{lemma}

These two lemmas and mean value theorem for $\beta(x,a)$ totally
dispense with the computation for \EqRef{sss}, and ease the computation
for \EqRef{A.1.5}, as follows:
\begin{lemma}\label{lem:up}
\[
\EqRef{sss}= \max\{  K_{a_c(p)}(q) \ \colon K_{a_c(p)}(q)<0,\ (p,q)\in\ad^2, \ p\prec q\}=K_{a_c(\tilde p)}(\tilde q)\change{\le\gap}{}.\]
\end{lemma}
  \begin{proof}  In the argument of $\max$, we have only to consider
   $\prec$-maximal $p$ such that $p\prec q$, by
   Lemma~\ref{lem:monotonic}. Since $q\in\ad$, there {is} a strictly increasing sequence
   $r_1<\cdots<r_k$ $({1\le k\le3})$ of integers in an interval $[3,\ {\gonality}]$
   and positive integers $n_1,\ldots,n_k$ such that $q$ consists of
   $n_i$ number of $r_i$ $(1\le i\le k)$. Then $\tilde p\in\ad$
   is obtained from $\tilde q\in\ad$ by replacing  exactly one component
   $j$ of $\tilde q$  with  $j-1$, because $\tilde p$ is maximal.
   Thus, \change{}{by Lemma~\ref{aux},}
   {$K_{a_c(\tilde p)}(\tilde q)=   \beta(j-1,a_c(\tilde p)) -
   \beta(j,a_c(\tilde p))$, which is $ -
   \frac{\partial\beta}{\partial x}(\xi, a_c(\tilde p))$ for some $\xi$
   $(j-1<\xi<j)$}\change{}{\ by mean value theorem for $\beta$}. By
   Lemma~\ref{lem:concave} and Lemma~\ref{lem:monotonic}~\EqRef{assert:ac monotone},
   $j=\gonality$ and $a_c(\tilde p)$ is the maximum among such $\tilde
   p$.
   This completes
   the proof.
  \end{proof}

 \begin{lemma}\label{lem:cncv} For $p=(f_1,f_2,f_3),\ q=(g_1,g_2,g_3)\in\ad$,
$$K_{a_c(f_1,f_2,f_3)}(g_1,g_2,g_3)<0\implies cncv(f_1,f_2,f_3,g_1,g_2,g_3).$$
\end{lemma}
\begin{proof} By Lemma~\ref{aux} and $a_c(p)>0$, we only need to verify the following:
   \begin{eqnarray}\label{lem:c}
    \left\{
   \begin{array}{l}
 \{l,m,n\}=\{1,2,3\},\ m<n,\ f_l\ge g_l\ \mbox{and}\
  f_m-g_m\ge g_n-f_n>0\\
\implies\sum_{i=1}^3 \left(\beta(f_i,a) - \beta(g_i,a)\right)>0\ (a>0).\end{array}\right.
\end{eqnarray}
We prove \EqRef{lem:c} as follows: {Suppose $l=1$. Then,
  by}  Lemma~\ref{lem:concave}~\EqRef{assert:increasing and concave in
 x}, the premise $f_l\ge g_l$ {of \EqRef{lem:c}} implies
\begin{equation}
 \label{u1}\sum_{i=1}^3 \left(\beta(f_i,a) - \beta(g_i,a)\right)\ge\sum_{i=2}^3(\beta(f_i,a)- \beta(g_i,a)).
\end{equation} By mean value theorem,
\begin{equation}
 \label{u2}\sum_{i=2}^3(\beta(f_i,a)-
  \beta(g_i,a))=\frac{\partial\beta}{\partial
  x}(\xi_2,\,a)(f_2-g_2)-\frac{\partial\beta}{\partial x}(\xi_3,\,
  a)(g_3-f_3)
\end{equation} for some $\xi_2$
  $(g_2<\xi_2<f_2)$ and $\xi_3$
  $(f_3<\xi_3<g_3)$, where the {premise} {of \EqRef{lem:c}} implies $f_2>g_2$ and $g_3>f_3$.
  By $p,q\in\ad$,
  $f_2\le f_3$. Thus,
  $\xi_2<\xi_3$. By Lemma~\ref{lem:concave}~\EqRef{assert:increasing and
  concave in x}, $\frac{\partial\beta}{\partial x}(\xi_2,\, a)>\frac{\partial\beta}{\partial x}( \xi_3,\, a)\change{}{>0}$. By the assumption
  $f_2-g_2\change{>}{\ge}g_3-f_3>0$, we have $\EqRef{u2}>0$. Hence, by
 \EqRef{u1}, we have \EqRef{lem:c}. {The cases $l\ne1$ are
 proved similarly.}
\end{proof}

\medskip
To sum up, the algorithm for $\kappa$ saves the computation time as follows: to find $(p,q)$ that attains
$\kappa$, the algorithm first sets $(p,q):=(\tilde p,\, \tilde q)$ and
then checks all the $(p,q)\in\ad^2$ satisfying the four conditions in
\EqRef{A.1.5} with updating $(p,q)$ whenever $0>K_{a_c(p)}(q)$ is
greater than the tentative maximum $K_{a_c(p)}(q)$.  When the computed
value of $K_{a_c(p)}(q)$ is less than the tentative maximum, the
algorithm dispenses computing $K_{a_c(p')}(q)$ for $p'\prec
p$\change{}{, because of
Lemma~\ref{lem:monotonic}~\EqRef{assert:curvature defined} and
\EqRef{assert:ac monotone}}.

Next, in order to eradicate the second difficulty, i.e., the round-off
difficulty,
\begin{itemize}
 \item
we will consider $\vartheta$ instead of
$K_{a_c(p)}(q)=2\pi - \vartheta$, and

 \item we will compute the fraction in the argument of $\arccosh$ in
       \EqRef{def:alpha}, based on Kahan's
	 method~\cite{kahan14:_miscal_area_angles_needl_trian,Goldberg91}
	 for Heron's formula for Euclidean triangles. Kahan's method is robust against loss of
	 significant digits according to the numerical experiment
	 \cite[Table~1]{kahan14:_miscal_area_angles_needl_trian}.
\end{itemize}

Now, our algorithm for $\gap$ consists of `Initialization and
subroutines' and `Main routine'.

\def\enddo{\textbf{end do}}
\def\endif{\textbf{end if}}
\def\For{\textbf{for}}
\def\If{\textbf{if}}
\def\Do{\textbf{do}}
\def\Else{\textbf{else}}
\def\Return{\textbf{return}}
\def\then{\textbf{then}}
\def\uIf{\textbf{if}}
 \begin{algorithm}
  \SetKwInput{Initialization}{Initialization and subroutines}
    \LinesNumbered
    \Initialization{}
For $3\le i\le \gonality, 3\le j<k\le    \gonality$, allocate memory for $C_i, S_i, D_{j,k}$, compute the values and let $C_i:=\cos(\pi/i), \quad   S_i:=C_i^2, \quad    D_{j,k}:=C_{k}-C_{j} $ ;
$\alpha(i,j,k)$:   {\Return} $\displaystyle \frac{8S_{i}S_{j}S_{k}}{(C_{k}+(C_{j}+C_{i}))    (C_{i}-D_{j,k}) (C_{i}+D_{j,k})    (C_{k}+D_{i,j})}$ ;
$\vartheta(A,i,j,k)$:        

{\Return} $\displaystyle \left(\arccos\left(\change{\frac{4     S_{i}}{A} - 1}{1-\frac{4     S_{i}}{A} }\right)+ \arccos\left(\change{\frac{4 S_{j}}{A} - 1}{1-\frac{4     S_j}{A}}\right)\right)  + \arccos\left(\change{\frac{4    S_{k}}{A} - 1}{1-\frac{4    S_{k}}{A}}\right)$ ; \end{algorithm}
By
	 \cite[Section~4]{kahan14:_miscal_area_angles_needl_trian} and
	 \EqRef{withintwice}, at most one bit of $D_{f_1,f_2}$ and
       	 $D_{f_2,f_3}$ is rounded off. Because of the first line of the
	 algorithm, the algorithm computes
	  $\cos(\pi/f_i)$, $\cos^2(\pi/f_i)$, and
	 $\cos(\pi/f_k)-\cos(\pi/f_j)$ at most once.

\medskip
 The subroutines $\alpha$ and $\vartheta$ of the algorithms have
the following properties:
    \begin{lemma}\label{lem:computation}
     For
 $p,q\in\ad$ and $A\ge2$,
  \begin{enumerate}
\item \label{real} The subroutine
   $\vartheta(A,p)$ returns a real number.

   \item  \label{acpaux}
	  $\alpha(p)=\cosh a_c(p)+1>2$.

   \item \label{vartheta:acpaux}
$ K_{\arccosh(\alpha(p)-1)}(q)=2\pi - \vartheta(\alpha(p),q)$.
  \end{enumerate}
\end{lemma}

 \begin{proof} 
  Let $p=(f_1,f_2,f_3)\in\ad$.
   By  $A\ge2$,
  $-1 \le 4\cos^2(\pi/f_i)/A - 1 \le 2\cos^2(\pi/f_i) - 1 =
  \cos(2\pi/f_i) \le 1.$ Thus, $\arccos(4\cos^2(\pi/f_i)/A -1)$ is a
  real number. Therefore, \EqRef{real} holds.
\EqRef{acpaux}
is due to Lemma~\ref{lem:monotonic}~\EqRef{assert:ac} and $a_c(p)>0$.
  \EqRef{vartheta:acpaux}. Let $q=(g_1,g_2,g_3)\in\ad$ and
let $\varphi$ be $\arcsin(\cos(\pi/g_i)/\cosh(a_c(p)/2))>0$. Then
  $\cos2\varphi=1-2\sin^2\varphi =
  1-2(\cos(\pi/g_i)/\cosh(a_c(p)/2))^2=1-2\cos^2(\pi/g_i)/((\cosh
  a_c(p)+1)/2)$. By \EqRef{acpaux} of this Lemma,
  $\cos2\varphi=1-4\cos^2(\pi/g_i)/\alpha(p)=1-4\cos^2(\pi/f).$ Therefore,
\[
2\arcsin\frac{\cos(\pi/g_i)}{\cosh(a_c(p)/2)}
 =\arccos\left(
\change{4\frac{\cos^2(\pi/g_i)}{\aux(p)} -
1}{1-4\frac{\cos^2(\pi/g_i)}{\aux(p)}} \right).\]
  The conclusion follows from Definition~\ref{def:monotonic}.
  This completes the proof.
 \end{proof}

\medskip  In computing the maximum $\kappa$, the main routine runs over many
of $(p,q)\in\ad^2$, so it uses the following function to run over all $q=(g_1,g_2,g_3)$ of $\ad$. Recall
Proposition~\ref{Hig}.
\begin{df}\label{def:dtwo and dthr}
 For $i$~$(3\le i \le\gonality)$, let $$m_2(i):= 7\ (i=3);\ 5\ (i=4,5);\ 6\ (i=6); \ i \ (\mbox{otherwise}).$$

For $i,j$~$(3\le i \le \gonality$, $m_2(i)\le j \le \gonality)$,  let
 $$m_3(i,j) :=
\left\{\begin{array}{cccc}
 43\ (i = 3,j = 7);& 25\ (i =  3,j = 8); & 19\ (i = 3,j = 9); \\
16\ (i = 3,j = 10); & 14\ (i = 3,j = 11);& 13\ (i = 3,j =  12); \\
21\ (i = 4,j = 5); & 13\ (i = 4,j = 6);& 10\ (i = 4,j = 7);\\
 9\ (i = 4,j = 8);   & 11\ (i = 5,j = 5); & 8\ (i = 5,j = 6);\\
 7\ (i = 5,j =7);  &   7\ (i = 6,j = 6);    & j\
 \left(-\frac{1}{2}+\frac{1}{i}+\frac{2}{j}<0\right);\\
\change{j}{\infty}\ \left(\mbox{otherwise}\right).
       \end{array}\right.
 $$
\end{df}

We set $\min\emptyset=\infty.$ 
\begin{lemma}\label{lem:AdP3} Let $i,j,k$ be integers. Then, $(i,j,k)\in\ad$
    if and only if
    \[
          3\le i\le \gonality, \quad
    m_2(i)\le j\le \gonality,\quad \mbox{and}\quad m_3(i,j)\le k\le \gonality.
    \]
   \end{lemma}
   \begin{proof}
Proposition~\ref{Hig} implies that
$m_2(i)$ is $\min\{j\colon \exists k.\,(i,j,k)\in \ad\}$, and that
    $m_3(i,j)$ is $\min\{k\colon (i,j,k)\in \ad\}$. 
\end{proof}

The main routine to compute the constant $\kappa$ is the following:

\begin{algorithm}\setcounter{AlgoLine}{4}\SetKwInput{Initialization}{Main routine}
    \LinesNumbered
    \Initialization{}
 {$\tilde{P}:=(\gonality-1,{\gonality},{\gonality})$ ;                \qquad\qquad   // cf. Lemma \ref{lem:up}
 $\tilde{Q}:=({\gonality},{\gonality},{\gonality})$ ;                  \qquad\qquad\qquad// cf. Lemma \ref{lem:up}
 $T:=\vartheta(\alpha(\tilde{P}), \tilde{Q})$ ;                 \qquad\qquad\qquad // Tentative minimum of $\vartheta$
 \For{ $g_1=3,\ldots, \gonality$} \Do
  \For{ $g_2=m_2(g_1),\ldots, \gonality$} \Do                        \qquad  // cf. Lemma \ref{lem:AdP3}
    \For{ $g_3=m_3(g_1,g_2),\ldots, \gonality$} \Do                  \qquad  // cf. Lemma \ref{lem:AdP3}
      $B$\_is\_computed:=\false ;     	// $B=\alpha(g_1,g_2,g_3)$ is not computed yet
      \For{$f_1=g_3-1,g_3-2,\ldots,3$\label{GrandFor}} \Do  // cf. remark just after the algorithm
        \For{ $f_2=\gonality,\gonality-1,\ldots, m_2(f_1)$} \Do
          $m$:=$\max$($g_1+1,$ $m_3(f_1,f_2))$ ;                 // \EqRef{sss} requires $p\not\prec q$
          \For{ $f_3=\gonality,\gonality-1,\ldots,m$} \Do
            \If { $\exists i,j. (f_i < g_i$ \& $f_j>g_j)$ \qquad         // \EqRef{sss} requires $p\not\prec q$
             \& $cncv(f_1,f_2,f_3,g_1,g_2,g_3)$ \label{cnd}}  \qquad\qquad\quad\qquad// cf. Lemma \ref{lem:cncv}
            \then {
               \If {$B$\_is\_computed=\false} \then
                  $B:=\alpha(g_1,g_2,g_3)$;                           //  $\cosh(a_c(q))+1$
                  $B$\_is\_computed:=\true
               \endif
              $A:=\alpha(f_1,f_2,f_3)$ ;                              //  $\cosh(a_c(p))+1$
                \If {$A<B$} \then {                 // \EqRef{sss} requires $a_c(p)< a_c(q)$
                      $k:=\vartheta(A,g_1,g_2,g_3)$ ;
                      \uIf {$k<T$} \then {
                          $T:=k$; $\tilde{P}:=(f_1,f_2,f_3)$ ; $\tilde{Q}:=(g_1,g_2,g_3)$ ;}
 }}                    \endif // \If $k<T$
                   \endif // \If $A<B$
                \endif //  \If $\exists \cdots$
              \enddo // \For $f_3 \cdots$
            \enddo // \For $f_2 \cdots$
          \enddo // \For $f_1 \cdots$
        \enddo // \For $g_3 \cdots$
      \enddo // \For $g_2 \cdots$
    \enddo // \For $g_1 \cdots$
    \Return $2\pi - T,\tilde{P},\tilde{Q}$\qquad\qquad    // $\gap,p,q$ such that $K_{a_c(p)}(q)=\gap$}
\end{algorithm}

\medskip
 The for-loop of $f_1$
 of the main routine decreases $f_1$  by $-1$, and thus does not consider case $g_3\le f_1\le
 \gonality$. However, this algorithm correctly computes the maximum $\kappa$ of
$ K_{a_c(p)}(q)<0$ such that \emph{both} of $ p$ and $q$ are in $\ad$
 and $a_c(p)<a_c(q)$, because of the following:
If $f_1\ge g_3$, then
	$(f_1,f_2,f_3)\succeq (g_1,g_2,g_3)$. Thus
	$a_c(f_1,f_2,f_3)\ge a_c(g_1,g_2,g_3)$ by
	Lemma~\ref{lem:monotonic}~\EqRef{assert:ac monotone}. Because
 $\kappa$ is greater than or equal to
  \EqRef{sss}, the algorithm does not have to check the case $g_3\le f_1\le\B{2}{3}$.


\begin{thebibliography}{10}
\expandafter\ifx\csname url\endcsname\relax
  \def\url#1{\texttt{#1}}\fi
\expandafter\ifx\csname urlprefix\endcsname\relax\def\urlprefix{URL }\fi
\expandafter\ifx\csname href\endcsname\relax
  \def\href#1#2{#2} \def\path#1{#1}\fi

\bibitem{Nevanlinna70}
R.~Nevanlinna, Analytic functions, Translated from the second German edition by
  Phillip Emig. Die Grundlehren der mathematischen Wissenschaften, Band 162,
  Springer-Verlag, New York-Berlin, 1970.

\bibitem{Stone76}
D.~Stone, {A combinatorial analogue of a theorem of Myers}, Illinois J. Math.
  20~(1) (1976) 12--21.

\bibitem{Gromov87}
M.~Gromov, Hyperbolic groups, in: Essays in group theory, no.~8 in Math.
  Sci. Res. Inst. Publ., Springer, New York, 1987, pp. 75--263.

\bibitem{Ishida90}
M.~Ishida, Pseudo-curvature of a graph, in: lecture at Workshop on topological
  graph theory, Yokohama National University, 1990.

\bibitem{Zuk97}
A.~{\.{Z}}uk, {On the norms of the random walks on planar graphs}, Annales de
  l'institut Fourier 47~(5) (1997) 1463--1490.

\bibitem{Woess98}
W.~Woess, A note on tilings and strong isoperimetric inequality, Math. Proc.
  Camb. Phil. Soc. 124 (1998) 385--393.

\bibitem{Higuchi01}
Y.~Higuchi, Combinatorial curvature for planar graphs, J. Graph Theory 38~(4)
  (2001) 220--229.

\bibitem{BP01}
O.~Baues, N.~Peyerimhoff, Curvature and geometry of tessellating plane graphs,
  Discrete Comput. Geom. 25~(1) (2001) 141--159.

\bibitem{HJL02}
O.~H{\"a}ggstr{\"o}m, J.~Jonasson, R.~Lyons, {Explicit Isoperimetric Constants
  and Phase Transitions in the Random-Cluster Model}, The Annals of Probability
  30~(1) (2002) 443--473.

\bibitem{LPZ02}
S.~Lawrencenko, M.~D. Plummer, X.~Zha, {Isoperimetric Constants of Infinite
  Plane Graphs}, Discrete {\&} Computational Geometry 28~(3) (2002) 313--330.

\bibitem{HigS03}
Y.~Higuchi, T.~Shirai, {Isoperimetric Constants of $(d,f)$-Regular Planar
  Graphs}, Interdisciplinary Information Sciences 9~(2) (2003) 221--228.

\bibitem{SunYu04}
L.~Sun, X.~Yu, {Positively curved cubic plane graphs are finite}, Journal of
  Graph Theory 47~(4) (2004) 241--274.

\bibitem{RBK05}
T.~R{\'e}ti, E.~Bitay, Z.~Kosztol{\'a}nyi, {On the polyhedral graphs with
  positive combinatorial curvature}, Acta Polytechnica Hungarica 2~(2) (2005)
  19--37.

\bibitem{BP06}
O.~Baues, N.~Peyerimhoff, Geodesics in non-positively curved plane
  tessellations, Adv. Geom. 6~(2) (2006) 243--263.

\bibitem{DeVosMohar07}
M.~DeVos, B.~Mohar, {An analogue of the Descartes-Euler formula for infinite
  graphs and Higuchi's conjecture}, Tran. Amer. Math. Soc. 7 (2007) 3287--3300.

\bibitem{ChenChen08}
B.~Chen, G.~Chen, {Gauss-Bonnet Formula, Finiteness Condition, and
  Characterizations of Graphs Embedded in Surfaces}, Graphs and Combinatorics
  24~(3) (2008) 159--183.

\bibitem{Zhang08}
L.~Zhang, {A result on combinatorial curvature for embedded graphs on a
  surface}, Discrete Mathematics 308~(24) (2008) 6588--6595.

\bibitem{Chenbf09}
B.~Chen, {The Gauss-Bonnet formula of polytopal manifolds and the
  characterization of embedded graphs with nonnegative curvature}, Proc. Amer.
  Math. Soc. 137~(5) (2009) 1601--1611.

\bibitem{Keller10}
M.~Keller, The essential spectrum of the {L}aplacian on rapidly branching
  tessellations, Mathematische Annalen 346 (2010) 51--66.

\bibitem{KP11}
M.~Keller, N.~Peyerimhoff, Cheeger constants, growth and spectrum of locally
  tessellating planar graphs, Math. Z. 268~(3-4) (2011) 871--886.

\bibitem{Keller11}
M.~Keller, {Curvature, geometry and spectral properties of planar graphs},
  Discrete \& Computational Geometry 46~(3) (2011) 500--525.

\bibitem{Oh17}
B.-G. Oh, {On the number of vertices of positively curved planar graphs},
  Discrete Mathematics 340 (2017) 1300--1310.

\bibitem{Gh17}
L.~Ghidelli, On the largest planar graphs with everywhere positive
  combinatorial curvature, arXiv:1708.08502 (2017).

\bibitem{HJL15}
B.~Hua, J.~Jost, S.~Liu, {Geometric analysis aspects of infinite semiplanar
  graphs with nonnegative curvature}, Journal f{\"u}r die reine und angewandte
  Mathematik 700 (2015) 1--36.

\bibitem{BuragoBuragoIvanov01}
D.~Burago, Y.~Burago, S.~Ivanov, A course in metric geometry, no.~33 in
  Graduate Studies in Mathematics, American Mathematical Society, Providence,
  RI, 2001.

\bibitem{MR2127379}
A.~D. Alexandrov, Convex polyhedra, Springer Monographs in Mathematics,
  Springer-Verlag, Berlin, 2005, translated from the 1950 Russian edition by N.
  S. Dairbekov, S. S. Kutateladze and A. B. Sossinsky, With comments and
  bibliography by V. A. Zalgaller and appendices by L. A. Shor and Yu. A.
  Volkov.

\bibitem{Thurston76}
W.~Thurston, Geometry and topology of 3-manifolds, Princeton lecture notes,
  1976.

\bibitem{ChowLuo03}
B.~Chow, F.~Luo, {Combinatorial Ricci flows on surfaces}, J. Differential Geom.
  63 (2003) 97--129.

\bibitem{LuoGuDai08}
F.~Luo, X.~Gu, J.~Dai, Variational principles for discrete surfaces, Vol.~2 of
  Advanced Lectures in Mathematics, International Press, Somerville, MA; Higher
  Education Press, Beijing, 2008.

\bibitem{GeXu16DGA}
H.~Ge, X.~Xu, {2-dimensional combinatorial Calabi flow in hyperbolic background
  geometry}, Differ. Geom. Appl. 47 (2016) 86--98.

\bibitem{GeXu17IMRN}
H.~Ge, X.~Xu, {A discrete Ricci flow on surfaces with hyperbolic background
  geometry}, Int. Math. Res. Not. 2017~(11) (2017) 3510--3527.

\bibitem{AkamaHuaSu18}
Y.~Akama, B.~Hua, Y.~Su, Areas of spherical polyhedral surfaces with regular
  faces, arXiv:1804.11033 (2018).

\bibitem{Bolza1887}
O.~Bolza, On binary sextics with linear transformations into themselves, Amer.
  J. Math. 10~(1) (1887) 47--70.

\bibitem{KKSV16}
K.~Katz, M.~Katz, M.~M. Schein, U.~Vishne, Bolza quaternion order and
  asymptotics of systoles along congruence subgroups, Experimental Mathematics
  25~(4) (2016) 399--415.

\bibitem{KS06}
M.~G. Katz, S.~Sabourau, An optimal systolic inequality for {CAT{$(0)$}}
  metrics in genus two, Pacific J. Math. 227~(1) (2006) 95--107.

\bibitem{iordanov17:_implem_delaun_trian_bolza_surfac}
I.~Iordanov, M.~Teillaud, Implementing {D}elaunay triangulations of the {B}olza
  surface, in: B.~Aronov, M.~J. Katz (Eds.), 33rd International Symposium on
  Computational Geometry (SoCG 2017), Vol.~44 of Leibniz International
  Proceedings in Informatics (LIPIcs), Schloss Dagstuhl--Leibniz-Zentrum fuer
  Informatik, Dagstuhl, Germany, 2017, pp. 44:1--44:15, article No.~44.

\bibitem{MR2836763}
R.~Nicholson, J.~Sneddon, New graphs with thinly spread positive combinatorial
  curvature, New Zealand J. Math. 41 (2011) 39--43.

\bibitem{Ol17}
P.~R. Oldridge., Characterizing the polyhedral graphs with positive
  combinatorial curvature, Master's thesis, Department of Computer Science,
  University of Victoria, https://dspace.library.uvic.ca/handle/1828/8030 (May
  2017).

\bibitem{Knuth:1997:ACP:270146}
D.~E. Knuth, The Art of Computer Programming, Volume 2 (3rd Ed.): Seminumerical
  Algorithms, Addison-Wesley Longman Publishing Co., Inc., Boston, MA, USA,
  1997.

\bibitem{kahan14:_miscal_area_angles_needl_trian}
W.~Kahan, Miscalculating area and angles of a needle-like triangle, Available
  via http://http.cs.berkeley.edu/\verb+~+wkahan/Triangle.pdf (Sep. 2014).

\bibitem{Goldberg91}
D.~Goldberg, What every computer scientist should know about floating-point
  arithmetic, ACM Comput. Surv. 23~(1) (1991) 5--48.

\end{thebibliography}
\end{document}